\documentclass{amsart}

\usepackage{amsmath}

\newtheorem{theorem}{Theorem}[section]
\newtheorem{lemma}[theorem]{Lemma}

\theoremstyle{definition}

\theoremstyle{remark}
\newtheorem{remark}[theorem]{Remark}

\numberwithin{equation}{section}

\begin{document}


\title [ ]{Determinantal approach to a proof of the Riemann hypothesis}


\author[]{}
\author[]{John Nuttall}
\address{}
\curraddr{}

\thanks{The assistance provided by the multiple precision arithmetic
  package written by David H. Bailey is much appreciated.}




\begin{abstract}
We discuss the application of the determinantal method to the proof of the Riemann Hypothesis. We start from the fact that, if a certain doubly infinite set of determinants are all positive, then the 
hypothesis is true.  This approach extends the work of Csordas, Norfolk and Varga in 1986, and makes extensive use of the results described by Karlin. We have discovered and proved or conjectured relations in five areas of the problem, as summarized in the Introduction. Further effort could well lead to more progress.
\end{abstract}

\maketitle


\newtheorem{prop}{Proposition}

\section{INTRODUCTION} {\label{secA}

\subsection{}

It is now over 150 years ago that the Riemann Hypothesis (RH) was first stated. Its proof or disproof is generally thought to be a very difficult problem, since many eminent mathematicians have so far been unable to solve it.

In this document\footnote{April 2011}    we\footnote{email: jnuttall@uwo.ca} suggest that the determinantal approach described below  allows the RH problem to be broken down into many perhaps simpler problems, which might facilitate a soloution  of the whole. The determinantal approach aims at the proof of the inequalites of (\ref{21f}). The method is an extension of the work in 1986 of 
 Csordas, Norfolk and Varga  ~\cite{CNV}  (see also Csordas and Varga  ~\cite{CV}), who proved (\ref{21f}) for the case of order $r=2$. We are not aware of any other significant progress on the determinantal approach since then.

In many instances we use the computer to explore the problem for possible relationships, to uncover apparent behavior of the functions involved as their variables change, and to carry out algebraic operations that would otherwise be extremely tedious. Thus we arrive at two types of result. In some cases, often stated in lemmas, we have rigorous statements, although often they were  discovered by calculation. Alternatively sometimes  we can only conjecture a general result, normally as an extension of patterns uncovered by calculation. In three cases we present explicit conjectures for important proposed results, although there are others that are implicit.

After a survey of previous relevant work in Sec.\:2, we report new developments in five topics related to the determinantal method. In Sec.\:3 we point out that the method used to prove the case of order $r=2$ has a natural generalization to higher order provided by the sign-regularity ideas described by Karlin ~\cite{K}, most importantly the appropriate constant sign of the Wronskian of $\Phi(u)$ of order $r$. Here $\Phi(u)$, defined in (\ref{21a}), is intimately connected with the Riemann $\zeta$-function. We refer to our recent reports ~\cite{JN2}, ~\cite{JN3} that use two different methods that prove the validity of  (\ref{21f}) for $r=3$ and $n\ge2$.  Numerical studies make it almost certain that the same holds for order $r=4$, except that the range of validity is then $n\ge3$. These studies also make it certain that the generalized method does not work for $r\ge5$. Of course it can be shown by rigorous numerical calculations that (\ref{21f}) does hold for the excluded cases for $r=3,4$, and indeed we have not found an example where the condiition fails in more than 100,000 tests - otherwise this report would have no point.

In Sec.\:4 we propose an extension of the generalized method used in Sec.\:3. The key element in that method is a kernel $K(u,v)=\Phi(u+v)$.  In Conjecture  \ref{cnj1} and the associated discussion we propose that the cumulants $\{\Psi_m(u)\}$ of $\Phi(u)$ defined in  (\ref{41a}), if used in the method in place of $\Phi(u)$, will lead to the proof of (\ref{21f}) for each $r$, provided that $n>\eta(r)$, where  $\eta(r)$ is a non-decreasing function of $r$. This conjecture is supported by numerical calculations presented in Sec.\:4.4. Thus, even if the Conjecture were true, for all $r>2$ there would be a set of values of $n$, increasing in number with $r$, for which the method would be unable to prove (\ref{21f}).

In Sec.\:5 we discuss one such case, $r=3,\:\:n=2$, and show that $D(2,3)>0$ will hold provided that $q(u,v)<0,\:\:u,v>0$, where $q(u,v)$ is a fairly simple function of $u,v$, analogous to the Wronskian, that is defined in   (\ref{54b}). The condition on $q(u,v)$ remains to be proved. If that can be accomplished, this example would provide hope that analogous methods will handle the other exceptional cases.

In Sec.\:6 we turn to an examination of the Wronskian of $\Phi(u)$ defined in  (\ref{61a}). The series  (\ref{21a}) defining $\Phi(u)$ converges rapidly, the more so as $u$ increases. Thus we can obtain a good approximation to the Wronskian for large $u$ if we replace the series by its first term. That leads to the study of the determinant $W_r(u)$ given by  (\ref{61c}) - examples up to $r=7$ are listed in  (\ref{61d}). From its definition it is clear that $W_r(u)$ is a polynomial of order $r^2$, but the examples suggest that the lowest $r(r-1)/2$ coefficients are zero, and likewise for the highest $r(r-1)/2$ coefficients. The elements in  (\ref{61c}) are polynomials $\{p_k(y)\}$ given by a recurrence relation  ~\cite[p.\,184]{CV}. We have found two alternative representations for the coefficients of the polynomials  $\{p_k(y)\}$. Using the first or lower representation we have shown in Lemma\:6.1 that the behavior observed for the  lowest $r(r-1)/2$ coefficients in $W_r(u)$ holds for all $r$. In Conjecture 2 we propose that the observed behavior of the highest $r(r-1)/2$ coefficients in $W_r(u)$ also applies for all $r$. We expect that the second or upper representation will be useful in the proof of this Conjecture. 

It is also clear from Sec.\:2.8 and ~\cite{JN3} that the behavior of the Wronskian (\ref{61a})  for $r=2,3$ is dominated by the first term of  (\ref{21a}). In the case of $r=3$ the contribution from the second term in the series is significant. This contribution may also be expressed in terms of determinants that are polynomials ~\cite[(3.6), p.\:6]{JN3}. Thus we see that the behavior of 
 $W_r(u)$  and a few corresponding  later determinants is crucial to the proof that $D(n,r)>0$ for $r=2,3,4$  (except for the exceptional cases). It might well be that the above representations are a key factor that will lead to a deeper understanding of why the RH holds.

In Sec.\:7 we discuss the question of the behavior of $D(n,r)$ as $n\rightarrow{\infty}$ for $n$ fixed. Our original approach ~\cite{JN4}, ~\cite{JN5}  was to use the Laplace method to approximate the integrals appearing in the definition of $D(n,r)$. We found that a proof that $D(n,r)>0$ for large $n$ is probably feasible provided that the truth of Conjecture 3 can be demonstrated. This conjecture is an algebraic relation not special to the RH. In Sec.\:7.3 we suggest an alternative approach that is also probably feasible. It is based on the properties of the determinant  $W_r(u)$, in particular on the requirement that the highest non-zero coefficient of 
$W_r(u)$ has an appropriate sign $\epsilon_r=(-1)^{r(r-1)/2}$ as stated in Conjecture\:2.

\section{PREVIOUS RESULTS}  {\label{secB}

\subsection{}

The foundations of the determinantal approach may be found in several sources, perhaps most importantly the book of Karlin ~\cite{K}. A key function  $\Phi(u)$, in the notation of  ~\cite{CNV} and  ~\cite{CV}, is defined by
\begin{equation}\label{21a}		  
\Phi(u)=\sum_{m=1}^{\infty}(2{m^4}{\pi^2}{e^{9u}}-3{m^2}{\pi}e^{5u})\exp{(-m^2}{\pi}e^{4u}).
\end{equation}
The function $\Phi(u)$ is even and analytic for all $u$ with ${\mid}{\Im{u}}{\mid}<{\pi}/4$ (see ~\cite{CNV}.

Let the Fourier cosine transform of $\Phi(u)$ be 
\begin{equation}\label{21b}   
\Xi(t)  =\int_{0}^{\infty\!\!\!}{du}\Phi(u)\cos(ut).
\end{equation}
It is well known that the RH is equivalent to the statement that all the zeros of $\Xi(t)$ are real. 

In 1925 P\'{o}lya therefore  raised the question of what properties of the  function $\Phi(u)$ are sufficient to lead to a proof that all the zeros of $\Xi(t)$ are real, but he presented no answers. To investigate this question we  set $z=-t^2$ and $F(z)=\Xi(t)$, so that equation (\ref{21b}) leads to
\begin{equation}\label{21c}   			  
\Xi(t)=F(z)  =\sum_{n=0}^{\infty}\beta_{n}z^n.
\end{equation}
Here the series coefficients  are
\begin{equation} \label{21d}   
\beta_{n} =
\frac{1}{\Gamma\left(2n+1\right)}\int_{0}^{\infty}\!\!{du}\Phi(u){u^{2n}},\:\:\:n=0,1,\ldots,.
\end{equation}

We define the matrix
\begin{equation}\label{21e}		
 B(i,j)
  =
  \left\{
  \begin{array}{ll}
    \beta_{j-i} , & j\ge i ; 
    \\
    0, &  j<i ;
    \\
   \end{array}
  \right. i,j=0,\, 1,\, 2,\, \ldots.
\end{equation}
As Karlin ~\cite[p.\:393, p.\:412]{K} explains, the RH requirement that all the zeros of $F(z)$ are real and negative is equivalent to the conditions
\begin{equation} \label{21f}		
D(n,r)>0,\:\:\:\: n=0,1,\ldots;\:\:r=1,2...,
\end{equation}
where we denote the minors $D(n,r)$ of order $r$ by
\begin{equation} \label{21g}				
D(n,r)=\det[B(i,j+n)]_{i,j=1,\ldots,r}.
\end{equation}
The proof of the inequalities  (\ref{21f})   is the set of simpler problems referred to above, and we call this approach the determinantal method.

\subsection{}

In the case of $r=1$ in  (\ref{21f}) the result is trivially correct, as all coefficients  (\ref{21d}) are positive, since $\Phi(u)>0$ for all real $u$.

For $r=2$ we have for all $n\ge0$
\begin{equation}\label{22a}     
D(n,2)=\left|
\begin{array}{cc}
\beta_n &\beta_{n+1} \\
\beta_{n-1} &\beta_n \\
\end{array} \right|.
\end{equation}
so that the condition (\ref{21f}) means that
\begin{equation}\label{22b}			
\beta_n^2>\beta_{n-1}\beta_{n+1},\:\:\:n=0,1,2,\ldots
\end{equation}
The proof of ({\ref{22b}}) may be regarded as the first step in the program we are suggesting.

In 1986 Csordas, Norfolk and Varga  ~\cite{CV} (see also Csordas and Varga [2])  proved the validity of the Tur\'{a}n inequalities (postulated about 50 years earlier by P\'{o}lya), which may be written as
\begin{equation} \label{22c}				
\beta_n^2>\frac{(n+1)}{n}\beta_{n-1}\beta_{n+1},\:\:\:n=1,2,\ldots
\end{equation}
Since $\frac{(n+1)}{n}>1$ it  follows that (\ref{22b})  is correct for $n=1,2,\ldots$, and the case $n=0$ is obvious. Thus the proof of the first problem  ({\ref{22b}}) is implicit in the results of Csordas, Norfolk and Varga.

\subsection{}

We now outline a proof of (\ref{22b}) using the basic ideas of ~\cite{CNV} , but which is simpler than their method since the result is weaker, but adequate for our purposes. We define the kernel
\begin{equation}\label{23a}		
K(u,v)=\Phi(u+v),\:\:0\le{u,v}
\end{equation}
Suppose that 
\begin{equation}\label{23b}		
\lambda(t)=\int_{0}^{\infty}\!\!dv\phi(v,t)\Phi(v),
\end{equation}
with the kernel $\phi(v,t)$ being
\begin{equation}\label{23c}			
\phi(v,t)=\frac{v^{t-1}}{\Gamma(t)},\:\:\:t>0.
\end{equation}
Note from (\ref{21d})  that
\begin{equation}\label{23d}			
\beta_n=\lambda(2n+1)\:\:\:n\ge0.
\end{equation}

Let the kernel $\Lambda(s,t)$ be defined by
\begin{equation}\label{23e}			
\Lambda(s,t)=\lambda(s+t),\:\:\:s,t>0.
\end{equation}
The argument underlying the development of   ~\cite[p.\,140]{K} shows that, for appropriate values of $s, t$, 
\begin{equation}\label{23f}						
\int_{0}^{\infty}\:du\int_{0}^{\infty}\:dv\phi(u,s)\Phi(u+v)\phi(v,t)
=\int_{0}^{\infty}dv\phi(v,s+t)\Phi(v)
\end{equation}
Thus we obtain
\begin{equation}\label{23g}					
\Lambda(s,t)=
\int_{0}^\infty{du}\int_{0}^\infty{dv}\phi(u,s)K(u,v)\phi(v,t).
\end{equation}

\subsection{}

We make considerable use of the concepts of sign-regularity and the compound kernel. Following Karlin  ~\cite{K} suppose that we have a kernel $Q(x,y)$. Let $X$ be a linearly ordered set, discrete or continuous. For a given positive integer $p$ the open simplex $\Delta_p\left(X\right)$ is
\begin{equation}\label{24a}			
\Delta_{p}\left(X\right)  =\left\{ \underline{x}=\left(x_{1},x_{2},\ldots,x_{p}\right)| 
           x_1<x_2<\ldots<x_p:x_i\in{X}        \right\},   
\end{equation}
and similarly for $Y$, $\underline{y}$.

As in  ~\cite[p.\,12]{K} the compound kernel $Q_{\left[{p}\right]}\left(\underline{x},\underline{y}\right)$ is defined by
\begin{equation}\label{24b}					
Q_{\left[{p}\right]}\left(\underline{x},\underline{y}\right)=\left| \begin{array}{cccc}
Q\left(x_1,y_1\right) &Q\left(x_1,y_2\right) &\ldots&Q\left(x_1,y_p\right) \\
Q\left(x_2,y_1\right) &Q\left(x_2,y_2\right) &\ldots&Q\left(x_2,y_p\right) \\
\vdots&\vdots&&\vdots\\
Q\left(x_p,y_1\right) &Q\left(x_p,y_2\right) &\ldots&Q\left(x_p,y_p\right) \\
\end{array} \right|.
\end{equation}
Suppose that the sequence $\epsilon_p=\left(-1\right)^{p(p-1)/2},\:\:p=1,\ldots,{r}$.  Then we say that $Q(x,y)$ is sign-reverse regular of order $r$ (i.e. $RR_r$) if $\epsilon_p   Q_{\left[{p}\right]}\left(\underline{x},\underline{y}\right)$
 is a non-negative function on  $\Delta_p\left(X\right) \times{\Delta_p\left(Y\right)}$ for each $p=1,2,\ldots,r$. Similarly, if $\epsilon_p=1,\:\:p=1,\ldots,{r}$, then $Q(x,y)$ is totally positive (i.e. $TP_r$). If the subscript $r$ is omitted in the designation we mean that the sign-regularity applies to all values of $r$.

From two applications of the  basic composition formula (BCF) ~\cite[(2.5), p.\,17]{K} to equation  (\ref{23g}) we find that for any $p>0$
\begin{equation} \label{24c}				
\Lambda_{[p]}(\underline{s},\underline{t})=
\int_{0}^\infty\!\!\!{d\underline{u}}\int_{0}^\infty\!\!\!{d\underline{v}}\:
\phi_{[p]}  (\underline{u},\underline{s})K_{[p]}(\underline{u},\underline{v})\phi_{[p]}(\underline{v},\underline{t}),
\end{equation}
where \underline{u}, \underline{v}, \underline{s}, \underline{t} are defined as in (\ref{24a}).

\subsection{}

The relevance of the above relations becomes clear when we combine the information in (\ref{24c}), (\ref{23e}) and  (\ref{23d}). Suppose in (\ref{24c}) we choose $p=2$, and \underline{s}, \underline{t} such that
\begin{equation}\label{25a}			
s_{1}=t_{1}=n-\frac{1}{2},\:\:\:s_{2}=t_{2}=n+\frac{3}{2},\:\:\:\:n=1,2\ldots
\end{equation}
The elements of the  determinant $\Omega(n)=\Lambda_{[2]}\left(\underline{s}, \underline{t} \right)$ are then
\begin{equation}\label{25b}		
\Omega(n)=\left[ \begin{array}{cc}
	\lambda(s_1+t_1) & \lambda(s_1+t_2)  \\
       \lambda(s_2+t_1)  &\lambda(s_2+t_2) \\
\end{array}\right],\:\:\: n=1,2,\ldots.
\end{equation}
Thus, using  (\ref{23d}), we have			
\begin{equation}\label{25c}		
\Omega(n)=\left| \begin{array}{cc}
	\beta_ {n-1} & \beta_{n} \\
        \beta_{n} & \beta_{n+1}\\
\end{array}\right|,\:\:\: n=1,2,\ldots.
\end{equation}

\subsection{}

It may be shown that the kernel $\phi(v,t)$ is $TP$, i.e. totally positive, which means that the determinant $\phi_{[p]}(\underline{u},\underline{s})$ is positive for $(\underline{u},\underline{s})$  in the appropriate simplex product  for any $p=1,2,\ldots$, and in particular for $p=2$.

 Now suppose that $K_{[2]}(\underline{u},\underline{v})<0$ for all values of the arguments (i.e. $K(u,v)$ is $RR_2$), then it follows from  (\ref{24c})  that
$\Lambda(u,v)$ is $RR_2$, so that $\Omega(n)<0,\:\:\: n=1,2,\ldots.$ We observe that, for $n\ge1$, the determinant $D(n,2)=-\Omega(n)$, since one determinant is the other with rows reversed.

The conclusion is that $D(n,2)>0,\:\:\:\: n=0,1,\ldots$, as required by the RH, provided that 
$K(u,v)=\Phi(u+v)$ is $RR_2.$

\subsection{}

An efficient method of testing for sign-regularity in the present case is described by Karlin  ~\cite{K}. We have 


\begin{lemma}   \label{L2a}			
Suppose that $\psi(x)$ is analytic in a neighborhood of $X=Y=(0,\infty)$, and that the kernel $K(x,y)=\psi(x+y),\:$  with $\:\:x,y\in (0,\infty)$. Define 
$w_p(u)=\det\left|  \psi^{(i+j-2)}(u)    \right|_{i,j=1}^p$.
If $\:\:{\epsilon_p} w_p(u)>0,\:\:u\ge0, \:\:\:p=1,2,\ldots,r\:\:\: then\: K(x,y) \:is \:RR_r$.
\end{lemma}
\begin{proof}    
This result is a special case of ~\cite[Theorem 2.6, p.\,55]{K}. The analyticity of $K(x,y)$ ensures that the differentiability requirements of the theorem are satisfied. The relation 
\begin{equation}\label{27b}			
\det\left|\frac{\partial^{i+j-2}}{\partial{x^{i-1}}\partial{y^{j-1}}}K(x,y)\right|_{i,j=1}^p=
\det\left|  \psi^{(i+j-2)}(u)    \right|_{i,j=1}^p,
\:\:\:u=x+y, 
\end{equation}
together with~\cite[(1.3), p.\,48]{K}  demonstrates that the requirements on the compound kernel appearing in the statement of the theorem hold.
\end{proof}

\newtheorem{rem2}{Remark}[section] 
\begin{rem2}
In the case of $r=2$ studied by ~\cite{CNV} and ~\cite{CV}, an alternative terminology is used. Instead of stating that $\psi(x+y)$ is $RR_2$ they require that $\log\psi(x)$ be strictly concave on $(0,\infty)$, i.e. that the second derivative of $\log{\psi(x)}$ be negative.  The two terms both lead to the condition (\ref{27c}).
\end{rem2}

To apply this lemma to the kernel $K(u,v)$ of (\ref{23a}) we set $\psi(x)=\Phi(x)$ and choose $r=2$. The required analyticity follows from  ~\cite[Theorem A, p.523]{CNV}.

For $p=1$ the condition of (\ref{27b}) requires that $\psi(u)>0,\:\:u\ge0$, since $\epsilon_1=1$ . This is true since $\Phi(u)>0,\:\:u\ge0$.

For $p=2$ the condition becomes
\begin{equation}\label{27c}			
-\left| \begin{array}{cc}
\Phi(u)	 & \Phi^{(1)}(u)\\
\Phi^{(1)}(u)	 &\Phi^{(2)}(u)\\      
\end{array}\right|=\Phi^{(1)}(u)^2-\Phi(u)\Phi^{(2)}(u)>0,\:\:u\ge0
\end{equation}
since  $\epsilon_1=-1$.

In their discussion of the  Tur\'{a}n inequalities, Csordas and Varga [2] study the case when $\psi(x)=\Phi(\surd{x})$. After "a series of involved, but elementary, estimates"  they prove that
\begin{equation}\label{27d}			
u[\Phi^{(1)}(u)^2-\Phi(u)\Phi^{(2)}(u)]+\Phi(u)\Phi^{(1)}(u)>0,\:\:u\ge0
\end{equation}
so that
\begin{equation}\label{27e}			
\Phi^{(1)}(u)^2-\Phi(u)\Phi^{(2)}(u)>-u^{-1}\Phi(u)\Phi^{(1)}(u)>0,\:\:u>0
\end{equation}
Inequality (\ref{27c}) follows since  $\Phi^{(1)}(u)<0,\:\:u>0$ and $\Phi^{(1)}(0)=0$ (see ~\cite{WN}).

Consequently we see that the work of ~\cite{CV} shows that $K(u,v)=\Phi(u+v)$ is $RR_2.$

\subsection{}

Here we present a self-contained, simpler proof of (\ref{27c}).  Using the notation of  ~\cite[p.\,184]{CV}  we have
\begin{equation}\label{28a}			
\Phi^{(j)}(u)=\sum_{n=1}^{\infty}a_n^{(j)}(u),\:\:\:j=0,1,\ldots
\end{equation}
where
\begin{equation}\label{28b}			
a_n^{(j)}(u)={\pi}{n^2}p_{j+1}({\pi}{n^2}{e^{4u}})\exp(5u-{\pi}{n^2}{e^{4u}}).
\end{equation}
With $y=\pi{e^{4u}}$ define 
\begin{equation}\label{28c}			
\Omega_j(y)=\Phi^{(j)}(u)\exp[-5u+{\pi}{e^{4u}}]/{\pi},\>\:\:\:j=0,1,2,
\end{equation}
so that (\ref{27c}) is equivalent to 
\begin{equation}\label{28d}			
W(y)=\Omega_0(y)\Omega_2(y)-\Omega_1(y)^2<0,\:\:\:y\ge\pi.
\end{equation}
We may write
\begin{equation}\label{28e}			
\Omega_j(y)=p_{j+1}(y)+\Upsilon_j(y),\:\:\:j=0,1,2.
\end{equation}
where, as in ~\cite{CV}
\begin{equation}\label{28f}	
\begin{array}{ccl}			
p_1(y )  & = & -3+2y\\
p_2(y)  & = & -15+30y-8y^2\\
p_3(y) & = & -75+330y-224y^2+32y^3.\\
\end{array}
\end{equation}
In (\ref{28e}) we define	
\begin{equation}\label{28g}			
\Upsilon_j(y)= {e^y}\sum_{n=2}^{\infty}{ n^2}p_{j+1}({n^2}y)\exp(-n^{2}y),\:\:\:j=0,1,2.
\end{equation}
Thus, substituting  (\ref{28e}) into (\ref{28d}), we obtain
\begin{equation}\label{28h}			
W(y)=W_2(y)+\sum_{i=1}^{5}T_i(y)
\end{equation}
where
\begin{equation}\label{28i}	
\begin{array}{ccl}			
W_2(y )  & = & p_1(y)p_3(y)-p_2(y)^2\\
T_1(y) & = & p_1(y)\Upsilon_2(y) \\
T_2(y) & = & p_3(y)\Upsilon_0(y) \\
T_3(y) & = & -2p_2(y)\Upsilon_1(y) \\
T_4(y) & = & \Upsilon_0(y)\Upsilon_2(y) \\
T_5(y) & = &- \Upsilon_1(y)^2  \\
\end{array}
\end{equation}

The proof of (\ref{27c}) is based on the following three Lemmas.


\begin{lemma}   \label{L2b}			
The function $W_2(y)<-843$ for $y\ge{\pi}$.
\end{lemma}
\begin{proof}    
Using the definitions (\ref{28f}) it may be shown that
\begin{equation}\label{28k}			
W_2(y)=16y(-15+12y-4y^2)=16y[-4(y-3/2)^2-6]
\end{equation}
Thus $W_2(y)$  decreases steadily for $y>1.5$, so that $W_2(y)<W_2(\pi) =-843.4...$ if $y>\pi$ and the lemma is proved.	
\end{proof}


\begin{lemma}   \label{L2c}		
For $y\ge{\pi}$ we have
$|p_k(y)|<2^{(2k-1)}y^{k},\:\:\:k=1,2,3$.
\end{lemma}
\begin{proof}    
For $k=1$ the proof is obvious. For $k=2$  the zero of 
\begin{equation}\label{28m}
p_2(y)+8y^2=30y-15
\end{equation}
is at $y=0.5$ Thus the function $30y-15$ is positive for $y\ge\pi>0.5$, so that 
\begin{equation}\label{28n}
p_2(y)>-8y^2,\:\:\:y\ge\pi
\end{equation}
Also we find that the zeros of 
\begin{equation}\label{28o}
p_2(y)-8y^2=-16y^2+30y-15
\end{equation}
are both complex. Thus  the function $-8y^2+30y-15$ is negative for $y\ge\pi>0.5$, so that 
\begin{equation}\label{28p}
p_2(y)<8y^2,\:\:\:y\ge\pi,
\end{equation}
from which the proof follows.
For $k=3$ the proof is similar. We find that the real zeros of $p_3(y)-32y^3$ are $< 1.193$, while the real zeros of  $p_3(y)+32y^3$ are $< 0.275$, which leads to the desired result.
\end{proof}


\begin{lemma}   \label{L2d}			
For $y\ge{\pi}$ and $j=0,1,2$ we have $|\Upsilon_j(y)|<2^{(4j+5)}C_j(y)y^{(j+1)} e^{-3y}  $ where $C_j(y)=(1-e^{-2y}{2^{(2j+4)}/2})^{-1}$.
\end{lemma}
\begin{proof}    
Csordas and Varga ~\cite[Lemma 3.1, p.\,185]{CV} have shown us how to prove this result.
Using (\ref{28g}) and	 the bounds on $p_j(y)$ of Lemma 2.3,  we find that
\begin{equation}\label{28r}
|\Upsilon_j(y)|<2^{2j+1}y^{j+1}e^y\sum_{n=2}^{\infty}n^{2j+4}\exp(-n^2{y})
\end{equation}
We may write
\begin{equation}\label{28s}
n^{2j+4}\exp(-n^2{y})=\frac{1}{\exp[n(ny-(2j+4)\log{n}/n)]}
\end{equation}
It is not difficult to show that, if $n\ge2$ and $2j+4<\pi$, then
\begin{equation}\label{28t}
(ny-(2j+4)\log{n}/n)>\log{K(n)}
\end{equation}
where $K(n)=e^{2y}{2^{-(2j+4)/2}}$.Thus
\begin{equation}\label{28u}
\sum_{n=2}^{\infty}n^{2j+4}\exp(-n^2{y})<\sum_{n=2}^{\infty}K(n)^{-n}=K^{-2}C_j(y)
\end{equation}
and the result follows.
\end{proof}

Combining the above information we are led to


\begin{lemma}   \label{L2E}
\begin{equation}\label{28v}			
\left|
\begin{array}{cc}
\Phi(u)	 & \Phi^{(1)}(u)\\
\Phi^{(1)}(u)	 &\Phi^{(2)}(u)\\      
\end{array}\right|<0,\:\:u\ge0
\end{equation}
\end{lemma}
\begin{proof}    
Above we explained how (\ref{28d}) is equivalent to  (\ref{28v}).	To prove  (\ref{28d}) we use Lemmas \ref{L2b} - \ref{L2d} to provide upper bounds to the absolute values (where appropriate) of all the quantities listed in (\ref{28i}), which for $y\ge\pi$ leads to
\begin{equation}\label{28w}	
\begin{array}{cll}			
W_2(y )  & < -843.19 & =-843.19\\
|T_1(y)| & <2^{14}C_2(\pi){{\pi}^4}e^{-3\pi} &=132.76  \\
|T_2(y)| & <2^{10}C_0(\pi){{\pi}^4}e^{-3\pi} &=8.30 \\
|T_3(y)| & <2^{13}C_1(\pi){{\pi}^4}e^{-3\pi}  &=64.88  \\
|T_4(y)| & <2^{18}C_0(\pi)C_2(\pi){\pi}^4e^{-6\pi} &=0.17 \\
|T_5(y)| & < 0& =0\\
\end{array}
\end {equation}
From (\ref{28h}) we deduce that an upper bound for $W(y)$ is the sum of the bounds in (\ref{28w}), which is 	-635.80. Thus the lemma is proved.
\end{proof}

\newtheorem{rem5}{Remark}[section]
\begin{rem5}
As we find by doing precise numerical computations, this bound is not tight. A much closer bound could be obtained by extending the method to explicitly include  the terms corresponding to $n=2$ in (\ref{28e}). Our systematic approach may be well suited to extension to  orders $r=3,4$, as the large number of terms involved can be handled on a computer (but still preserving rigor).
\end{rem5}

\section{DOES THE OBVIOUS EXTENSION WORK?}  {\label{secC}

\subsection{}

The way in which we have described the work of  ~\cite{CNV}, ~\cite{CV}  in Sec.\:2 suggests an obvious method for attempting to extend the proof of the conditions  (\ref{21f})   to higher values of the order $r$. We can adjust the procedure of Sec.\:2.5, and Lemma \ref{L2a} applies to all values of $r$.    We need to determine whether the  kernel $K(u,v)$ of (\ref{23a}) has sign-regularity pattern   $RR_r$ for values of $r>2$.

An easy way to attack this question is via a numerical computation. With $\psi(u)=\Phi(u)$, the dervatives of $\psi(u)$ for use in (\ref{27b}) may be written as rapidly converging series obtained from (\ref{21a}), which can be approximated to an accuracy of hundreds of decimal places with little trouble. Calculations for a set of values of $u$ indicate that 
\begin{equation}\label{31a}
 \:\:{\epsilon_r} w_r(u)>0,\:\:u\ge0, \:\:\:r=1,2,\ldots,4
\end{equation}
but the inequality is definitely not true for order $r=5$.  If our code is not in error, then $K(u,v)$ is not $RR_5$, and the method of Sec.\:2 cannot be used to prove any of the inequalities of  (\ref{21f}) for $r\ge5$. 

We discuss a method that may to some extent overcome this difficulty in Sec.\:4, but next we report on the case of $r=3$.

\subsection{}

We have produced a proof ~\cite{JN2}) that
\begin{equation} \label{32a}		
D(n,3)>0,\:\:\:\: n=2,3,\ldots
\end{equation}
The method generally follows the pattern of ~\cite{CNV} and ~\cite{CV}, in that some extensive computations (not numerical, but computer assisted) are used to prove that the kernel $K(u,v)$ is $RR_3$. Equation  (\ref{24c}) with $p=3$ then shows that $\Lambda_{[3]}(\underline{s},\underline{t})$ is negative for all (\underline{s},\underline{t}) provided that the components of the arguments are all $>0$.

A more systematic method of proving the same result was described in a later report ~\cite{JN3}.

To  prove  (\ref{32a}) we  choose 
\begin{equation}\label{32b}
s_{1}=t_{1}=n-\frac{3}{2},\:\:\:s_{2}=t_{2}=n+\frac{1}{2},\:\:\:s_{3}=t_{3}=n+\frac{5}{2},\:\:\:\:n=2,3,\ldots
\end{equation}
and as in Sec. 2.5  use the relations  (\ref{24b}),  (\ref{23e}),  (\ref{23d}),  which give    
\begin{equation} \label{32c}
\Lambda_{\left[{3}\right]}\left(\underline{s},\underline{t}\right)
=    \left| \begin{array}{ccc}
        \beta_{n-2}   & \beta_{n-1} & \beta_{n}   \\
        \beta_{n-1}   & \beta_{n} & \beta_{n+1}   \\
        \beta_{n}   & \beta_{n+1} & \beta_{n+2}   \\	
\end{array} \right|<0,\:\:\:\:n=2,3,\ldots
\end{equation}
The restriction $n\ge2$ arises from the need to satisfy the condition $s>-1/2$   from (\ref{23c}).

Reversing the order of  the rows of (\ref{32c}), thereby changing its sign, produces the determinant $D(n,3)$ (remember that $\epsilon_3=-1$) , so that
\begin{equation} \label{32d}
D(n,3)=\left| \begin{array}{clcr}
        \beta_ {n} & \beta_{n+1} & \beta_{n+2} \\
        \beta_{n-1} & \beta_{n} & \beta_{n+1}	\\    
        \beta_{n-2} & \beta_{n-1} & \beta_{n}	\\
\end{array}\right|>0,\:\:\:\:n=2,3,\ldots
\end{equation}

Just as for order $r=2$ the case $n=0$ is trivial, but the case  $n=1$ is exceptional in that the  relation 
\begin{equation} \label{32e}
  D(1,3)=  \left| \begin{array}{ccc}
        \beta_ {1} & \beta_{2} & \beta_{3} \\
        \beta_{0} & \beta_{1} & \beta_{2}	\\    
        0 & \beta_{0} & \beta_{1}	\\
\end{array}\right|>0,
\end{equation}
does not follow from  the general method. Its validity may be checked by inserting the numerical values of $\beta_j$ listed by~\cite[p.\,540]{CNV}.

In Sec. 5 below we present some further observations on the case of $D(1,3)$.

\subsection{}

Suppose that our calculation leading to (\ref{31a}) is correct for $r=4$, which means that $K(u,v)$ is $RR_4$. Then we can apply reasoning similar to  that in Sec. 3.2. We choose
\begin{equation}\label{33a}
s_{j}=t_{j}=n-{5/2}+2j,\:\:\:\ j=0,1,2,3;  \:\:\:n=3,4\ldots
\end{equation}
with the result that
\begin{equation} \label{33b}
\Lambda_{\left[{4}\right]}\left(\underline{s},\underline{t}\right)
=    \left| \begin{array}{cccc}
       \beta_{n-3}   & \beta_{n-2} & \beta_{n-1}   & \beta_{n}  \\
        \beta_{n-2}   & \beta_{n-1} & \beta_{n}  & \beta_{n+1}   \\
        \beta_{n-1}   & \beta_{n} & \beta_{n+1}   & \beta_{n+2}  \\
        \beta_{n}   & \beta_{n+1} & \beta_{n+2}   & \beta_{n+3}  \\	
\end{array} \right|<0,\:\:\:\:n=3,4,\ldots
\end{equation}

In this case reversing the order of the rows does not change the sign of the dererminant, but, since $\epsilon_4=1$, we still obtain 
\begin{equation} \label{33c}
D(n,4)
=    \left| \begin{array}{cccc}
       \beta_{n}   & \beta_{n+1} & \beta_{n+2}   & \beta_{n+3}  \\	
       \beta_{n-1}   & \beta_{n} & \beta_{n+1}   & \beta_{n+2}  \\
        \beta_{n-2}   & \beta_{n-1} & \beta_{n}  & \beta_{n+1}   \\
        \beta_{n-3}   & \beta_{n-2} & \beta_{n-1}   & \beta_{n}  \\
 
\end{array} \right|>0,\:\:\:\:n=3,4,\ldots
\end{equation}

Now there are two non-trivial cases that are not covered by the general method, i.e. $D(1,4)$ and $D(2,4)$), but again a numerical calculation shows that these determinants are positive. We note that, in all the exceptional cases, one or more elements of the determinant is zero.

\section{CUMULANTS}{\label{secD}

\subsection{}

It appears that the answer to the question of the previous section is that the obvious generalization of the method of ~\cite{CNV} and    ~\cite{CV}  probably does allow a modest extension to orders 3 and 4. To proceed further we present a conjecture that could provide a partial proof of the RH inequalities for all orders. Even if the conjecture is correct, there will remain  to be studied an infinite set of exceptional cases such as occur for orders $3,4$. A beginning to this study is found in  Sec.\:\ref{secE}.

We introduce the cumulants $\{\Psi_m(u)\}$, where
\begin{equation}\label{41a}						
\Psi_m(u)=\int_{u}^{\infty}\!\!{du}\Psi_{m-1}\left(u\right),\:\:\:m=1,2,\ldots
\end{equation}
and
\begin{equation}\label{41b}					
\Psi_{0}(u)=\Phi(u).
\end{equation}	
There is a corresponding set of kernels
\begin{equation}\label{41c}
K\left(u,v;m\right)=\Psi_{m}\left(u+v\right),\:\:0\le{u,v};\:\:\:m=0,1,\ldots.
\end{equation}

As in (\ref{23f}) it may be shown that, for appropriate values of $m, s, t$, 
\begin{equation}\label{41d}				
\int_{0}^{\infty}\:du\int_{0}^{\infty}\:dv\phi(u,s)\Psi_m(u+v)\phi(v,t)
=\int_{0}^{\infty}dv\phi(v,s+t)\Psi_m(v)
\end{equation}
 Integrating by parts $m$ times leads to
\begin{equation}\label{41e}				
\int_{0}^{\infty}dv\phi(v,s+t)\Psi_m(v)
=\int_{0}^{\infty}dv\phi(v,s+t+m)\Phi(v)=\lambda(s+t+m)
\end{equation}
where $t, s>0$ and $m=0,1,2,\ldots$.
Thus we obtain
\begin{equation}\label{41f}				
\Lambda\left(s+m/2,t+m/2\right)=
\int_{0}^\infty{du}\int_{0}^\infty{dv}\phi\left(u,s\right)K\left(u,v;m\right)\phi\left(v,t\right).
\end{equation}

As for (\ref{24c}), using the compound kernel $K_{[p]}(\underline{u},\underline{v};m)$ corresponding to \:\:\:\:\:\: $K(v,u:m)$, we have
\begin{equation} \label{41g}					
\Lambda_{\left[{p}\right]}\left(\underline{s}+m\underline{w}/2  ,\underline{t}+m\underline{w}/2 \right)=
\int_{0}^\infty\!\!\!{d\underline{u}}\int_{0}^\infty\!\!\!{d\underline{v}}\:
\phi_{\left[p\right]}  \left(\underline{u},\underline{s}\right)K_{\left[p\right]}\left(\underline{u},\underline{v};m\right)\phi_{\left[p\right]}\left(\underline{v},\underline{t}\right),
\end{equation}
where \underline{u}, \underline{v}, \underline{s}, \underline{t} are defined as in  (\ref{24a})), and  \underline{w} is a vector with all components equal to unity.

\subsection{}

Now we adapt the procedure of  (\ref{32b}),  (\ref{33a}) to relate the compound kernel 
appearing in the LHS of (\ref{41g}) to some of the determinants needed in  (\ref{21f}). In (\ref{41g}) we choose $p=r>1$, and \underline{s}, \underline{t} such that
\begin{equation}\label{42a}				
s_j=t_j=\mu+2j,\:\:\:j=1,\ldots,r.
\end{equation}
The elements of the  determinant $\Lambda_{\left[{r}\right]}\left(\underline{s}+m\underline{w}/2  \underline{t}+m\underline{w}/2 \right)$ are then
\begin{equation}\label{42b}					
\lambda(2\mu+2i+2j+m),\:\:\:i,j=1,\ldots,r
\end{equation}
If the argument of $\lambda$ appearing in (\ref{42b}) is odd and positive for all entries, then we can use  (\ref{23e}) to express all elements of $\Lambda_{\left[{r}\right]}\left(\underline{s}+m\underline{w}/2  \underline{t}+m\underline{w}/2 \right)$ in terms of the coefficients $\beta_n$ given by (\ref{23d}). For use in  (\ref{41g}) the inequality  (\ref{23c}) requires that $\mu+2>0$. It follows that the determinant with elements given by (\ref{42b}) is the same as  $D(k,r)$ with the order of the rows reversed, and an appropriate choice of $k$. As before, the reversal of the rows changes the determinant by a factor of $\epsilon_r$.

It follows that all values of $k$ are possible so long as $k\ge{k_L}$. It is found that
\begin{equation}\label{42c}		
k_L=n_L+r-1,
\end{equation}
where  $n_L=m/2,\:\:m$ even:\:\:\:$(m+1)/2,\:\:m$ odd.

\subsection{}

Let us  assume the validity of the following conjecture    

\newtheorem{con1}{CONJECTURE}
\begin{con1}  \label{cnj1}
For any given order ${r}>1$, there is a lowest integer $m(r)$ such that   $K(u,v;m(r))$ is $RR_{r}$.
\end{con1}

\newtheorem{rem1}{Remark}[section] 
\begin{rem1}
Karlin  ~\cite[(1.17), p.\,102; Remark 5.6, p.\,128]{K}  shows that, because of the relation (\ref{41a}), if $K(u,v;m_1)$ is $RR_{r}$, then so is $K(u,v;m_2)$  for all $m_2>m_1$. Thus  the conjecture implies that, if $m\ge{m(r)}$, then $K(u,v;m)$ is $RR_{r}$. It also follows that, if $r_{2}>r_{1}$, then $m\left(r_{2}\right)\ge{m\left(r_{1}\right)}$.
\end{rem1}

The Conjecture means that the compound kernel $K_{\left[r\right]}\left(\underline{u},\underline{v};m(r)\right)$ has sign  $\epsilon_r$, so that, as before, the same applies to
$\Lambda_{\left[{r}\right]}\left(\underline{s}+m\underline{w}/2, \underline{t}+m\underline{w}/2 \right)$  from (\ref{41g}), since the other factors in the integrand are non-negative. 

The conclusion is that, given the Conjecture, we can prove the positivity of $D(n,r)$ needed in     (\ref{21f}) for all values of $n,r$ such that $n\ge{k_L}{(m(r))}$.  

\subsection{}

To examine in a non-rigorous numerical way whether the conjecture might be valid we have applied the method of Sec.\:3.1 to the kernels $K(x,y;m), \:\:m=0,1,2,\ldots;\:\:\:r=2,3,\ldots$, with each sequence going beyond $16$. The results are consistent with the Conjecture, some values of $m(r),\eta(r))$ being listed in the Table, where $\eta(r)=k_L(m(r))$ denotes the number of exceptional cases (not including $n=0$) for the indicated value of $r$.

\begin{equation}\label{44a}
\begin{array}{crrrrrrrrrrrrrrrrrrrr}
r &2 &3 & 4 & 5 &6 &7 &8 &9 &10 &11 &12 &13 &14 &15 &16& 17& 18& 19& 20 \\
m(r)&0 &0 & 0 & 1 &1 &1 &2 &4 &6 &7 &7 &7 &9 &11 &13 &15 &15&  15& 16 \nonumber\\
\eta(r)&1 &2 & 3&5 &6 &7 &8 &10 &12 &14 &15 &16 &18 &20 &22&24&25&26&27 \\
\end{array} 
\end{equation}

If the conjecture is correct, then the proof of the RH reduces to showing that
\begin{equation} \label{44b}						
D(n,r)>0,\:\:\:\: n=0,\ldots ,\eta(r);\:\:r=1,2...,
\end{equation}

\subsection{}

In some applications it might be helpful to have an explicit form for the cumulant  $\Psi_m(u)$ rather than the iterative relations (\ref{41a}).
Karlin  ~\cite{K} discussed cumulants (although apparently with a different aim in mind), and he showed that 
\begin{equation} \label{45a}
\Psi_m(u)=\frac{1}{\Gamma(m)}\int_{u}^\infty\!\!\!dt\Phi(t)(t-u)^{m-1},\:\:\:m\ge1.
\end{equation} 
This formula may be verified by differentiating with respect to $u$, which shows that
\begin{equation} \label{45b}
\Psi_m^{(1)}(u)=-\Psi_{m-1}(u),\:\:\:m\ge1,
\end{equation} 
consistent with  (\ref{41a}) and  (\ref{41b}).

\section{THE EXCEPTIONAL CASE $D(1,3)$} {\label{secE}

\subsection{}

In Sec.\:{\ref{secC} we explained that the relations (\ref{31a}) ${\epsilon_r} w_r(u)>0,\:\:u\ge0, \:\:\:r=1,2,3$ led to the result  (\ref{32a})  that $D(n,3)>0,\:\:\:\: n=2,3,\ldots$. In this section we describe an analogous technique, which shows that the relation $D(1,3)>0$ holds provided that $q(u,v)<0$ for all $u, v>0$, where $q(u,v)$ is defined in (\ref{54b}). This is the simplest example of an exceptional case.  

We modify the development of Sec.\,2.2 by replacing $\Phi(u)$ by $f(u)=\Phi^{(1)}(u)$, and defining the kernel $F(u,v)=f(u+v)$. Note that ~\cite[p.\,523]{CNV} states that $f(0)=0$ because $\Phi(u)$ is an even function of $u$. Thus, with $s, t>0$,    (\ref{23g}) becomes
\begin{equation}\label{51a}						
\int_{0}^{\infty}\!\!\!du\int_{0}^{\infty}\!\!\!dv\phi(u,s)f(u+v)\phi(v,t)
=\int_{0}^{\infty}dv\phi(v,s+t)f(v)=\omega(s+t),
\end{equation}
where
\begin{equation}\label{51b}						
\Omega(s,t)=\omega(s+t)=\int_{0}^{\infty}dv\phi(v,s+t)f(v).
\end{equation}

As in  (\ref{24c}) we apply the BCF with $p=3$, and  find that 
\begin{equation} \label{51c}						
\Omega_{\left[{3}\right]}\left(\underline{s},\underline{t}\right)=
\int_{0}^\infty\!\!\!{d\underline{u}}\int_{0}^\infty\!\!\!{d\underline{v}}\:
\phi_{\left[3\right]}  \left(\underline{u},\underline{s}\right)F_{\left[3\right]}\left(\underline{u},\underline{v}\right)\phi_{\left[3\right]}\left(\underline{v},\underline{t}\right),
\end{equation}
where \underline{u}, \underline{v}, \underline{s}, \underline{t}   are defined as before, including inequalities for \underline{u}, \underline{v} as in (\ref{24a}). We replace the RHS of  (\ref{51c}) using  ~\cite[(2.6)\:p.17]{K}, so that (\ref{51c}) becomes
\begin{equation} \label{51d}						
\Omega_{\left[{3}\right]}\left(\underline{s},\underline{t}\right)=
\int_{0}^\infty\!\!\!{d\underline{u}}\int_{0}^\infty\!\!\!{d\underline{v}}\:
\Theta(\underline{u},\underline{s}) F_{\left[3\right]}\left(\underline{u},\underline{v}\right)\Theta(\underline{v},\underline{t}),
\end{equation}
where
\begin{equation} \label{51e}						
\Theta(\underline{u},\underline{s})=\phi(u_1,s_1) \phi(u_2,s_2)\phi(u_3,s_3)
\end{equation}
and all the integrals now run over $(0,\infty)$.

\subsection{}

Now  in (\ref{51d}) we choose \underline{s}, \underline{t} with $z>0$ such that
\begin{equation}\label{52a}				
s_j=t_j=z+2j-2,\:\:\:j=1,\ldots,3.
\end{equation}
The elements of the  determinant $\Omega_{\left[{3}\right]}\left(\underline{s},  \underline{t} \right)$ are then
\begin{equation}\label{52b}				
\Delta(i,j)=\omega(2z+2i+2j-4),\:\:\:i,j=1,\ldots,3.
\end{equation}

The next step is to take the limit as $z\rightarrow{0}$ of both sides of (\ref{51c}). On the LHS consider first the  element 
\begin{equation}\label{52c}					
\Delta(1,1)=\omega(2z)=\int_{0}^{\infty}dv\frac{v^{2z-1}}{\Gamma(2z)} \Phi^{(1)}(v)
\end{equation}
Integrating by parts (differentiating $\Phi^{(1)}(v)$)  leads to
\begin{equation}\label{52d}						
\Delta(1,1)=\frac{-1}{\Gamma(2z+1)}\int_{0}^{\infty}dv{v^{2z}} \Phi^{(2)}(v)
\end{equation}
so that      
\begin{equation}\label{52e}					
\lim_{z \rightarrow 0} \Delta(1,1) = 0,
\end{equation}
since
\begin{equation}\label{52f}					
\int_{0}^{\infty}dv \Phi^{(2)}(v)=\Phi^{(1)}(\infty)-\Phi^{(1)}(0)=0
\end{equation}
on account of the properties of $\Phi^{(1)}(u)$.

For the remaining cases we have 
\begin{equation}\label{52g}					
\Delta(i,j)=\int_{0}^{\infty}dv\frac{v^{2z+2i+2j-5}}{\Gamma(2z+2i+2j-4)} \Phi^{(1)}(v),\:\:i,j=1,2,3;\:\:\:i+j>2.
\end{equation}
Integrating by parts (integrating $\Phi^{(1)}(v)$) and  setting $z=0$  leads to
\begin{equation} \label{52h}
\Delta(i,j)=\frac{-1}{\Gamma(2i+2j-3)}\int_{0}^{\infty}dv{v^{2i+2j-6} \Phi}(v)        
= \lambda(2i+2j-5) = -\beta_{i+j-3}
\end{equation}   


Thus for $z=0$ the elements of $\Omega_{\left[{3}\right]}\left(\underline{s},  \underline{t} \right)$ are the negative of the elements of the determinant $D(1,3)$ (see (\ref{32e})  )   with the order of the rows reversed.

\subsection{}

We now turn to the RHS of  (\ref{51d}). The function $1/\Gamma(z) =O(z)$ as $z\rightarrow{0}$ , so that in that limit the contribution to (\ref{51d}) from any region of integration with $u_1>u_0$ or $v_1>u_0$ for any positive $u_0$ will be zero. Since 
\begin{equation}\label{53a}					
\Gamma(z)=\int_{0}^{\infty}\!\!\!du\,u^{z-1}e^{-u},
\end{equation}
we end up with 


\begin{lemma}  \label{L5a}					
Apart from an inessential positive constant factor,
\begin{equation} \label{53c}					
      D(1,3)=\!\!
\int_{0}^\infty\!\!\!\!\!d{u_2}   \int_{0}^\infty\!\!\!\!\!d{u_3} \int_{0}^\infty\!\!\!\!\!d{v_2} \int_{0}^\infty\!\!\!\!d{v_3} 
\phi(u_2,2)\phi(u_3,4) F_{\left[3\right]}\left(\underline{u},\underline{v}\right)\phi(v_2,2)\phi(v_3,4)
\end{equation}
where $\underline{u}=(0,u_2,u_3),\underline{v}=(0,v_2,v_3)$
\end{lemma}
\begin{proof}
In case the above remarks are not considered to be an adequate proof of this lemma, we outline an alternative method. Expand each of the two determinants in  (\ref{53c}) into their four components (remember that the $(1,1)$ element of $F_{[3]}(\underline{u},\underline{v})=f(u_1+v_1)=\Phi^{(1)}(0)=0)$, and perform the integrals. It will be found that the four terms on each side of the equation are equal term by term. For example, consider the component from  
$F_{[3]}(\underline{u},\underline{v})$ of the form $f(u_2)f(v_2)f(u_3+v_3)$, which leads to
\begin{align}\label{53d}
&\int_{}^{}d{u_2}d{u_3}d{v_2}d{v_3}  
\phi(u_2,2)\phi(u_3,4) \phi(v_2,2)\phi(v_3,4)f(u_2)f(v_2)f(u_3+v_3)\\
&=C\left(\int_{}^{}d{u_2}u_2f(u_2)\right)\left(\int_{}^{}d{v_2}v_2f(v_2)\right)
\left(\int_{}^{}d{u_3}d{v_3}{u_3}^{3}{v_3}^{3}f(u_3+v_3)\right)\nonumber\\
&=C\left(\int_{0}^{\infty}duu\Phi^{(1)}(u)\right)^2
\left(\int_{0}^{\infty}du\int_{0}^{\infty}dvu^{3}{v}^{3}\Phi^{(1)}(u+v)\right)\nonumber
\end{align}\nonumber

The last line is easily shown to be a multiple of ${\beta_{0}^2}{\beta_{3}}$, where $\beta_n$ is given by  
(\ref{21d}). When the appropriate constant factors are included, a similar procedure for the other three terms leads to Lemma \ref{L5a}, after taking account of  (\ref{21g}).
\end{proof}

\subsection{}

To further examine the case of $D(1,3)$ define $Q(u_2,u_3,v_2,v_3)$ as
\begin{align}\label{54a}			
Q(u_2,u_3,v_2,v_3)&=F_{[3]}(\underline{u},\underline{v})
=\left| \begin{array}{ccc} 0 &f(v_2) &f(v_3) \\
f(u_2) &f(u_2+v_2) &f(u_2+v_3) \\
f(u_3) &f(u_3+v_2)&f(u_3+v_3) \\
\end{array} \right|
\end{align} 
where $\underline{u}=(0,u_2,u_3)$, $\underline{v}=(0,v_2,v_3)$, and $f(0)=0$ as mentioned above.  In analogy with the function $w_p(u)$ defined in Lemma 2.1 we introduce $q(u,v)$ as
\begin{align} \label{54b}	
q(u,v)&=\left| \begin{array}{ccc} 0 &f(v) &f^{((1)}(v) \\
f(u) &f(u+v) &f^{(1)}(u+v) \\
f^{(1)}(u) &f^{(1)}(u+v)&f^{(2)}(u+v), \\
\end{array} \right|
\end{align} 
which relates to the behavior of $Q(u_2,u_3,v_2,v_3)$  for ${u_2}\approx{u_3}$, 
 ${v_2}\approx{v_3}$.

In analogy with Lemma \ref{L2a} we have


\begin{lemma}   \label{L5b}			
Suppose that $q(u,v)<0$ for all $u, v>0$, and also that $f(u)<0$ for all $u>0$. Then $Q(u_2,u_3,v_2,v_3)>0$ for all $u_2>u_3$ and $v_2>v_3$.
\end{lemma}
\begin{proof}    
Make the definitions
\begin{equation} \label{54d}				
h(u,v)=\frac{f(u+v)}{f(u)}
\end{equation}
\begin{equation} \label{54e}			
R(u,v)=\frac{1}{f(v)^2}\left[f(v)\frac{\partial{h(u,v)}}{\partial{v}}-f^{(1)}(v)h(u,v)\right]
\end{equation}
\begin{equation} \label{54f}			
G(u_2,u_3,v)=\frac{1}{f(v)}\left[h(u_2,v)-h(u_3,v)\right]
\end{equation}

The following relations are easily checked.
\begin{equation} \label{54g}			
\frac{\partial}{\partial{v}}\left[\frac{1}{f(v)}h(u,v)\right]=\frac{1}{f(v)^2}R(u,v)
\end{equation}
\begin{equation} \label{54h}				
\frac{\partial{G(u_2,u_3,v)}}{\partial{v}}  = \frac{1}{f(v)^2}\left[R(u_2,v)-R(u_3,v)\right]
\end{equation}
\begin{equation} \label{54i}				
\frac{\partial{R(u,v)}}{\partial{u}}=\frac{-1}{f(u)^2}q(u,v)
\end{equation}
\begin{equation} \label{54j}				
Q(u_2,u_3,v_2,v_3)=f(u_2)f(u_3)f(v_2)f(v_3)\left[G(u_2,u_3,v_2)-G(u_2,u_3,v_3)\right]
\end{equation}

The assumption of Lemma \ref{L5b} together with (\ref{54i}) shows that
\begin{equation} \label{54k}		
\frac{\partial{R(u,v)}}{\partial{u}}>0
\end{equation}
 for all $u,v>0$. Thus (\ref{54h}) implies that 
\begin{equation} \label{54m}					
\frac{\partial{G(u_2,u_3,v)}}{\partial{v}} >0 
\end{equation}
for all $v>0$ if $u_2>u_3>0$.
This in turn means with (\ref{54j}) that 
\begin{equation} \label{54n}				
Q(u_2,u_3,v_2,v_3)>0
\end{equation}
if $u_2>u_3$ and $v_2>v_3$ as required in the statement of the lemma.
\end{proof}


\begin{lemma}   \label{L5c}
With the assumptions of Lemma \ref{L5b} then $D(1,3)>0$.
\end{lemma}
\begin{proof}    
From (\ref{54a}) it is clear that 
\begin{equation} \label{54q}
Q(u_3,u_2,v_2,v_3)=Q(u_2,u_3,v_3,v_2)=-Q(u_2,u_3,v_2,v_3).
\end{equation}
Thus, using  (\ref{54a}) for $F_{[3]}(\underline{u},\underline{v})$, with $C$ a positive constant, we can rewrite Lemma 3.1 as
\begin{equation} \label{54r}
      D(1,3)=C\!\!
\int_{0}^\infty\!\!\!\!\!\!d{u_2}   \int_{0}^{u_2}\!\!\!\!\!\!d{u_3} \int_{0}^\infty\!\!\!\!\!\!d{v_2} \int_{0}^{v_2}\!\!\!\!\!d{v_3}   \phi_{\left[2\right]}\left(\underline{u},\underline{s}\right)Q(u_3,u_2,v_2,v_3)\phi_{\left[2\right]}\left(\underline{v},\underline{t}\right),
\end{equation}
where $\underline{u}=(u_2,u_3)$,  $\underline{v}=(v_2,v_3)$ and $\underline{s}= \underline{t}=(2,4)$. For $u_2>u_3$ and  $v_2>v_3$ we have  
\begin{equation} \label{54s}
 \phi_{\left[2\right]}\left(\underline{u},\underline{s}\right)
 \phi_{\left[2\right]}\left(\underline{v},\underline{t}\right)
=\frac{u_2u_3v_2v_3}{6}
(u_3^{2}-u_2^{2})(v_3^{2}-v_2^{2})>0,
\end{equation}
so that (\ref{54r}) shows that $D(1,3)>0$ as required.
\end{proof}

\newtheorem{rem3}{Remark}[section] 
\begin{rem3}
The assumption in Lemma \ref{L5b} that $f(u)=\Phi^{(1)}(u)<0$ for $u>0$ may be disregarded. It was proved to be true by Wintner  ~\cite{WN}.
\end{rem3}


\section{COEFFICIENT REPRESENTATIONS} {\label{secF}


\subsection{}

We described in Sec.\:2.8 how the derivatives of $\Phi(u)$ can be written in terms of a set of polynomials  $p_k(y),\:\:k=1,2,\ldots,$ introduced by Csordas and Varga ~\cite[(3.9), p.\,185]{CV} . These CV polynomials play an important role in the study of the determinantal method. For instance consider the Wronskian
of general order $r$
\begin{equation}\label{61a}
w_{r}\left(u\right)= 
  \left| \begin{array}{cccc}
        \Phi(u)   & \Phi^{(1)}(u) & \ldots &  \Phi^{(r-1)}(u)  \\
        \Phi^{(1)}(u) & \Phi^{(2)}(u) &\ldots &  \Phi^{(r)}(u) \\
	\vdots & \vdots & & \vdots \\
        \Phi^{(r-1)}(u)  & \Phi^{(r)}ut) & \ldots & \Phi^{(2r-2)}(u)   	
\end{array} \right| 
\end{equation} 
Using the first term in the expansion of $\Phi^{(k)}$ we see that,  for large values of $u$,  we may approximate $w_{r}\left(u\right)$ by 
\begin{equation}\label{61b}
w_r(u){\approx}\left[\frac{exp[-5u+y]}{\pi}\right]^r{W_r(y)},
\end{equation} 
where
\begin{equation}\label{61c}
W_{r}\left(y\right)= 
  \left| \begin{array}{cccc}
        p_1(y)   & p_2(y) & \ldots &   p_r(y)  \\
        p_2(y) &  p_3(y) &\ldots &   p_{r+1}(y) \\
	\vdots & \vdots & & \vdots \\
        p_r(y)  &  p_{r+1}(y) & \ldots &  p_{2r-1}(y)  	
\end{array} \right| ,
\end{equation} 
with  $y=\pi{e^{4u}}$. We have seen in Sec.\:2.8 and ~\cite{JN3} how the properties of $W_r(y)$ were crucial in  proving  the positivity  of $D(n,r)$ for $r=2,3$. Later in Sec.\:7.3 we shall argue that the properties of  $W_r(y)$ could also form the basis of a proof that positivity applies for any fixed $r$ so long as $n$ is large enough.

The polynomials $p_j(y)$ of degree $j$ are defined by recurrence relations, starting with $p_1(y)=2y-3$, that are given in Csordas and Varga   ~\cite[p.\,184]{CV}. From the definition (\ref{61c}) it is clear  that $W_r(y)$ is a polynomial of degree at most $r^2$. 

 Using $\epsilon_r={(-1)}^{\mu(r)}$, where $\mu(r)=r(r-1)/2$, we find that
\begin{equation}\label{61d}	
\begin{array}{cl}		
{\epsilon_2}W_2(y)  & =240y-192y^2+64y^3  \\
{\epsilon_3}W_3(y)   & =  -860160y^3+737280y^4-294912y^5+65536y^6  \\
{\epsilon_4}W_4(y)   & =  190253629440{y^6}  -169114337280{y^7} +  72477573120{y^8} \\
 &      - 19327352832{y^9}   + 3221225472{y^{10}} \\
{\epsilon_5}W_5(y)   & =   -0.329167393077068 \:10^{19}{y^{10}}  +   0.299243084615516\:10^{19}{y^{11}}  \\
  &  -0.132996926495785\:10^{19}  {y^{12} } +  0.379991218559386\:10^{18}  {y^{13} }\\
  &   -0.759982437118771\:10^{17}  {y^{14} }+   0.101330991615836\:10^{17}  {y^{15} }\\
{\epsilon_6}W_6(y)   & =   0.538444964246560\:10^{28}  {y^{15} }       -0.497026120842978\:10^{28}  {y^{16} } \\
    & +  0.225920964019536\:10^{28} {y^{17} }        -0.669395448946772\:10^{27} {y^{18} }    \\
     &   + 0.143441881917165\:10^{27 } {y^{19} }     -0.229507011067465\:10^{26} {y^{20} }      \\
     & + 0.255007790074961\:10^{25} {y^{21} }   \\
 {\epsilon_7} W_7(y)   & =      -0.975629606681896\:10^{39}  {y^{21} }     +     0.910587632903103\:10^{39 }   {y^{22} }    \\        &-0.420271215186047\:10^{39  }    {y^{23} }    +     0.127354913692742\:10^{39 }  {y^{24} }         \\ &-0.283010919317204\:10^{38   } {y^{25} }     +       0.485161575972349\:10^{37  }   {y^{26} }      \\ &-0.646882101296465\:10^{36  }   {y^{27} }   +     0.616078191710919\:10^{35 }  {y^{28} }       \\
\end{array}
\end{equation}
Note that all the coefficients in (\ref{61c}) are integers, but we have displayed only the 15 leading digits.

The examples in (\ref{61d}) suggest three general  properties of the polynomials  $W_r(y)$.  We have proved one of these properties, which is stated in Lemma  \ref{L6a}. Some progress has been made on proving the other two properties, which are described in Conjecture  \ref{cnj2}. These properties all relate to the coefficients $\{\gamma(j,r)\}$ of the  polynomials $\{W_r(y)\}$ defined by
\begin{equation}\label{61e}				
{\epsilon_r}  W_r(y)=\sum_{j=0}^{r^2}\gamma(j,r)y^j.
\end{equation}
We have


\begin{lemma}   \label{L6a}
For $j=0,\ldots,\mu(r)-1$  the coefficients     $\gamma(j,r)=0,\:\: r\ge2$.
\end{lemma}
\begin{proof}  The proof, given in Sec.\:6.3, depends on  Lemma\:\ref{L6c} .
\end{proof}


\begin{con1}  \label{cnj2}
For a given $r>1$ the highest non-zero coefficient is $\gamma(\mu(r+1),r)$, and  $\gamma(\mu(r+1),r)>0$.
\end{con1}

\subsection{}

The above lemma depends on the existence of expressions for the coefficients of the CV polynomials ~\cite[p.\,184]{CV}, which are defined by
\begin{equation}\label{62a}				
p_k(y)=\sum_{j=0}^{k}d(j,k)y^j.
\end{equation}
There are two alternative representations of the coefficients $\{d(j,k)\}$, the first of which is described in this section.  We call it  the 'lower representation' because it is useful in proving Lemma  \ref{L6a}. We expect that the second, called the 'upper representation', will be useful in proving Conjecture \ref{cnj2}. It is described in Sec.\:6.4.


\begin{lemma}   \label{L6c}
For $k\ge1$  the coefficients     $d(j,k)$ are given by
\begin{equation}\label{62b}				
d(j,k)=\sum_{\eta=0}^{j}  c(j,\eta) (4\eta+5)^{k-1},\:\:\:j\le {k},
\end{equation}
where
\begin{equation}\label{62c}				
c(j,\eta)=\frac{{{(-1)}^{\eta+1}}(2\eta+3)}{{\eta}!(j-\eta)!},\:\:\:\eta\le{j}.
\end{equation}
Also, with $d(j,k)$ defined by (\ref{62c}),
\begin{equation}\label{62d}				
d(j,k)=0,\:\:\:j>k.
\end{equation}
\end{lemma}
\begin{proof}  The relation given by  ~\cite[p.\,184]{CV} is
\begin{equation}\label{62e}				
p_{k+1}(y)=4yp_k^{(1)}(y)+(5-4y)p_k(y), \:\:\:k=1,2,\ldots.
\end{equation}
Substituting (\ref{62a}) for $k\ge1$  leads to
\begin{equation}\label{62f}				
\begin{array}{cl}
d(0,k+1)&=5d(0,k)  \\
d(j,k+1)&=-4d(j-1,k)+[4j+5]d(j,k),\:\:j=1,\ldots,k \\
d(k+1,k+1)&=-4d(k,k) ,
\end{array}
\end{equation}
where the recurrence starts with
\begin{equation}\label{62g}				
d(0,1)=-3;\:\:\:d(1,1)=2.
\end{equation}

It is obvious that there is a unique solution for $d(j,k),\:\:j\le{k}$ obtained by iterating (\ref{62f}) on $k$ and starting from the conditions (\ref{62g}). Therefore the lemma will be proved if we can show that the definition (\ref{62b}) satisfies (\ref{62f}) and (\ref{62g}). 

For the first line of (\ref{62f}), the  formula  (\ref{62b}) shows that
\begin{equation}\label{62h}	
d(0,k)=-{3\times}{5^{k-1}}			
\end{equation}
so that 
\begin{equation}\label{62i}	
d(0,k+1)=5d(0,k) 			
\end{equation}
as required.

With a given value of $j$ the difference $J$  between the two sides of the second line of equation (\ref{62f}) is
\begin{equation}\label{62j}
J=\sum_{\eta=0}^{j} c(j,\eta)(4\eta+5)^{k-1}[4\eta+5-4j-5]  +
    4 \sum_{\eta=0}^{j-1} c(j-1,\eta)(4\eta+5)^{k-1}.				
\end{equation}
Since $(j-\eta)c(j,\eta)=c(j-1,\eta), \:\:\:\eta=0,\ldots,{j-1}$ it follows that $J=0$ as required.

To examine the third line we note that  the definition (\ref{62b}) leads  to
\begin{equation}\label{62k}	
\begin{array}{cl}
d(k,k)  &= \sum_{\eta=0}^{k}{ {(-1)}^{\eta+1}    (2\eta+3)(4\eta+5)^{k-1}}{  [{\eta!} {(k-\eta)!}]^{-1} }	\\
& ={2\times}{4}^{k-1} \sum_{\eta=0}^{k}{   {(-1)}^{\eta+1} [ {\eta}^k +P(\eta) ]  } \:  {  [{\eta!} {(k-\eta)!}]^{-1} },	
\end{array}
\end{equation}
where $P(\eta)$ is a polynomial in $\eta$ of degree $k-1$.
Now we define
\begin{equation}\label{62l}	
f(x)={k!} \sum_{\eta=0}^{k}   \frac{ x^{\eta}   }   {   {\eta!} {(k - \eta)!}     }=(1+x)^k.		\end{equation}
Therefore, with $j\le{k}$,
\begin{equation}\label{62m}	
\left(x\frac{d}{dx}\right)^j {f(x)} \left|_ {x=-1} \right. =   {k!} \sum_{\eta=0}^{k}       \frac{{\eta^j} (-1)^{\eta}   }   {   {\eta!} {(k - \eta)!}     }   ={(-1)^k}   { f^{ j)  }}  {(-1)}  
\end{equation}
since $f^{(j)}(-1)=0,\:\:\:j<k.$
It follows that
\begin{equation}\label{62n}	
 \sum_{\eta=0}^{k}  \frac{{\eta^k} (-1)^{\eta}   }   {   {\eta!} {(k - \eta)!}  }=
    (-1)^k,\:\:j=k;\:\:\:=0,\:\:j<k.
\end{equation}
 Using this result in  (\ref{62k}) shows that
\begin{equation}\label{62o}	
d(k,k)=(-1)^{k+1}{2\times}{4}^{k-1}, k=1,2,\ldots 
\end{equation}
since the contribution of the term containing $P(\eta)$ will be zero. Together the above results prove (\ref{62b}) of the lemma. Equation (\ref{62d}) follows in the same way using (\ref{62j}).
\end{proof}


Some values of $c(j,\eta)$ are shown in  (\ref{62p})
\begin{equation}\label{62p}
\begin{array}{rcccccl}
	j & \eta=0  & \eta=1 & \eta=2  & \eta=3 &  \eta=4 \\
	0 & -3\\
	1 & -3 & 5  \\ 
	2 & -3/2 & 5 & -7/2    \\
	3 & -3/6 & 5/2 & -7/2 & 9/6  \\
	4 & -3/24 & 5/6 & -7/4 & 9/6 & -11/24
\end{array}
\end{equation}

\subsection{}

We now use the representation of the coefficients $d(j,r)$ demonstrated in Lemma \ref{L6c} to prove Lemma  \ref{L6a}. Define a set of column vectors $\underline{v}(\eta)$ of dimension $r$ by
 $\underline{v}(\eta)=[v_1(\eta),v_2(\eta),\ldots,v_r(\eta)] ^T  $
where
\begin{equation}\label{63a}	
  v_{\alpha}(\eta)=(4\eta+5)^{\alpha-1},\:\:\:\alpha=1,\ldots,r.		
\end{equation}
Also let the column vector $\underline{V}(\beta,y)= [p_{\beta}(y),p_{\beta+1}(y),\ldots
p_{\beta+r-1}(y)] ^T,\:\:\:\beta=1,\ldots,r,\;\;\:$ so that
\begin{equation}\label{63b}	
	  W_r(y)=det\left[  \underline{V}(\beta,y)
		 \right]	_ {\beta=1,\ldots,r}
\end{equation}
Then we have, using (\ref{62a}) and  (\ref{62b}),
\begin{equation}\label{63c}
  \underline{V}(\beta,y) =   \sum_{j=0}^{r+\beta-1}\underline{H}(j,\beta){y^j},\:\:\:
\beta=1,2,\ldots,r, 
\end{equation}
where
\begin{equation}\label{63d}
  \underline{H}(j,\beta) =   \sum_{\eta=0}^{j}c(j,\eta)(4\eta+5)^{\beta-1}
 \underline{v}(\eta),	\:\:\:
\beta=1,2,\ldots,r. 
\end{equation}

To prove Lemma \ref{L6a} we note that we have shown in Lemma \ref{L6c} that each column $\underline{V}(\beta,y)$  of the determinant  $W_r$ is a polynomial in $y$ of degree $r+\beta-1$. Each coefficient in these polynomials is a linear combination of column vectors $\underline{v}(\eta)$. Note that $\eta{\le}j$ for the $\underline{v}(\eta)$ vectors appearing in the coefficient of $y^j$.

The determinant $W_r(y)$ may be expressed as a sum of the determinants of all combinations of $r$ column vectors  $\underline{v}(\eta)$ obtained by choosing one  $\underline{v}(\eta)$ from each column  $\underline{V}(\beta,y),\:\:\:\beta=1,\ldots,r$. In this process, we consider two $\underline{v}(\eta)$ with the same value for $\eta$ to be different if they appear in terms with different $j$.

 Every example of the resulting determinants will include as a factor $ det\left[  \underline{v}(\eta(\beta))  \right]	_ {\beta=1,\ldots,r}$, where $\eta(\beta)$ is the value of $\eta$ chosen for the term in column $\beta$. If for any two different values $\beta=\beta_1,
\beta_2$ it happens that $\eta(\beta_1)=\eta(\beta_2)$, then the determinant will be zero since it contains two identical columns  $\underline{v}(\eta(\beta_1))$ and
 $\underline{v}(\eta(\beta_2))$.

If no pair of  $\eta(\beta_1)$, $\eta(\beta_2)$ are equal then 
\begin{equation}\label{62r}	
 \sum_{\beta=1}^{r}  \eta(\beta){\ge}  \sum_{i=0}^{r-1}i=(r(r-1))/2=\mu(r)
\end{equation}
and conversely, if 
\begin{equation}\label{62s}	
 \sum_{\beta=1}^{r}  \eta(\beta)<\mu(r)
\end{equation}
then  must be such an equal pair.
If $j(\beta)$ is the value $j$ associated with $\eta(\beta)$ it follows, since   $j(\beta){\ge}
\eta(\beta)$, that there is an equal pair if 
\begin{equation}\label{62t}	
 \sum_{\beta=1}^{r}  j(\beta)<m,
\end{equation}
thus proving the claim of Lemma \ref{L6a}.


\begin{remark} 
Using similar arguments it may be proved that the coefficient $\gamma(\mu(r),r)$ defined in  (\ref{61e}) is given by
\begin{equation}\label{62pp}	
 \gamma(m,r)=\left\{ {\prod}_{j=0}^{r-1}c(j,j)    \right\}
\left[     det \left| (5+4i)^j\right|_{i,j=0,\ldots,r-1}  \right]^2
\end{equation}
This result may be extended to a few higher coefficients without much difficulty.
\end{remark}

\subsection{}

We now turn to the second of the alternative representations for the coefficients of the CV polynomials, the 'upper representation'. For $k=5,9$ we  define an integer function $s(i,j),\:\:i=1,2,\ldots$ by the initial conditions
\begin{equation}\label{64a}	
\begin{array}{rl}
 s(i,j)  & =0,\:\:i\le0,\:\:j=1,2,\ldots    \\
 s(i,j)  & =0,\:\:j<0,\:\:i=1,2,\ldots    \\
 s(i,0)  & =1,\:\:i=1,2,\ldots    \\

 s(1,j)& =k^j,\:\:j=1,2,\ldots     \\
 s(i,i)  & =2\times{(-4)^{i-1}},\:\:i=1,2,\ldots    \\
\end{array}
\end{equation}
and then recursively by
\begin{equation}\label{64aa}
 s(i+1,j) = s(i,j)+s(i+1,j-1)(k+4i),\:\: i=1,2,\ldots,\:\:j=1,2,\ldots
\end{equation}
We obtain some properties of $\{s(i,j)\}$ in


\begin{lemma}   \label{L6d}
For $i,j\ge{1}$ let $\underline{\nu}=(\nu_1,\nu_2,\ldots,\nu_i)$,   where $\{\nu_{\alpha}\}$ are integers such that $0\le\nu_{\alpha}\le{j}$ \: and \:  $ \sum_{\alpha=1}^{i}{\nu_\alpha}=j$.  For given $\{i,j\}$ let\: $\underline{\nu}(m),\:\:m=1,2,\ldots,N(i,j)$ be the list of all possible partitions with $i\ge1,\:\:j\ge0$, thus defining $N(i,j)$.  Then
\begin{equation}\label{64b}	
N(i,j)=\frac{(i+j-1)!}{j!{(i-1)!}}
\end{equation}
and
\begin{equation}\label{64c}	 s(i,j)=\sum_{m=1}^{N(i,j)}  \left(     \prod_{\alpha=1}^{i}                         (  k+(4({\alpha}-1)  )^  {  \nu_{\alpha}(m)  }  
  \right).
\end{equation}
\end{lemma}
\begin{proof}  
The proof is by induction. First we note that, for $j=1$, there are $i$  partitions  with
 ${\nu}_{\alpha}(m)=1,\:\:m=\alpha;\:\:=0,\:\:m\ne{\alpha}$ for $m,\alpha=1,2,\ldots,i$.
Thus $N(i,1)=i$, in agreement with (\ref{64b}), and  (\ref{64c}) becomes	
\begin{equation}\label{64d}	
 s(i,1)=\sum_{m=1}^{i} (k+4(m-1)),
\end{equation}
which leads to
\begin{equation}\label{64e}	
 s(i+1,1)- s(i,1)  =k+4i,\:\:i=1,2,\ldots
\end{equation}
This is consistent with  (\ref{64a}) for $j=1$, bearing in mind that $s(i,0)=1$.

Similarly, when $i=1$, there is just one partition, which $\nu_1(1)=j$. In this case   (\ref{64b}) gives $N(1,j)=1$, and  (\ref{64c}) leads to $s(1,j)=k^j$, both of which are consistent with  (\ref{64a}).

Next  we  investigate the sequence $\{s(i,2),\:\:i=2,3,\ldots \}$. Consider first $i=2$. From 
(\ref{64aa}) we see that $s(2,2)$ should be related to $s(1,2)$ and $s(2,1)$. We know that  
$N(1,2)=1$ and that the single partition for $s(1,2)$ is $\underline{\nu}(1)=(2)$, with $s(1,2)=k^2$. Also $N(2,1)=2$ and that the two partitions for $s(2,1)$ are $\underline{\nu}(1)=(1,0),\:\:\underline{\nu}(2)=(0,1)$, with $s(2,1)=k+[k+4]$.

According to the definition in the statement of the lemma, the partitions corresponding to $s(2,2)$ must have the form $\underline{\nu}(m)=(\nu_1(m),\nu_2(m))$ where 
\begin{equation}\label{64f}	
0\le\nu_{\alpha}(m)\le{2},\:\:\alpha=1,2;\:\::m=1,2,\ldots,N(2,2)
\end{equation}
and
\begin{equation}\label{64g}	
\nu_1(m)+\nu_2(m)=2,\:\:m=1,2,\ldots,N(2,2).
\end{equation}

This means that $N(2,2)=3$ and that the 3 partitions must be
\begin{equation}\label{64h}	
\underline{\nu}(1)=(2,0),\:\:\underline{\nu}(2)=(1,1),\:\:\underline{\nu}(3)=(0,2)
\end{equation}
 with 
\begin{equation}\label{64i}	
s(2,2)=[k^2]+[k(k+4)]+[(k+4)^2=[k^2] +[k+(k+4)](k+4)
\end{equation}

We see that (\ref{64i}) corresponds precisely to (\ref{64aa}) with $i=1,\:j=2$, which reads 
\begin{equation}\label{64j}	
s(2,2)=s(1,2)+s(2,1)(k+4)
\end{equation}

In general $s(i+1,j)$ is the sum of all monomials of degree $j$ constructed from products of numbers in the sequence $k,k+4,k+8, \ldots,k+4i$. The equation  (\ref{64aa}) declares that this sum is obtained by taking the sum all  the monomials in $s(i,j)$, which are also of degree $j$ but include only those formed from the sequence $k,k+4,k+8, \ldots,k+4(i-1)$. To that must be added the sum of all monomials of degree $j-1$  in $s(i+1,j-1)$ constructed from products of numbers in the sequence $k,k+4,k+8, \ldots,k+4i$, which are then mulitiplied by $ k+4i$. All the monomials in $s(i+1,j)$ will appear from one or other of these two terms, and there will be no repeats.

In a similar fashion we obtain 
\begin{equation}\label{64k}	
N(i+1,j)=N(i,j)+N(i+1,j-1)
\end{equation}
which is satisfied by (\ref{64b}), so that the Lemma is proved. 
\end{proof}
	
We are now in a position to state and prove the   CV polynomial coefficients 'upper representation'. Let $s_5(i,j),s_9(i,j)$ denote the function $s(i,j)$ of   (\ref{64a}), (\ref{64aa}), with the subscript being the value of $k$.


\begin{lemma}   \label{L6e}
The coefficients $\{d(j,k)\}$ defined in  (\ref{62a}) are given by
\begin{equation}\label{64l}	
	\begin{array}{cl}
d(i,i+j)   &=-3\times{(-4)^i}s_5(i+1,j-1)  \\
& +2\times{(-4)^{(i-1)}}s_9(i,j),\:\: i=0,1,\ldots;\:\:j=
0,1,\ldots
\end {array}
\end{equation}
\end{lemma}
\begin{proof}
For the special cases where $\{i=0,\:\:j=1,2,\ldots\}$ or $\{i=j,\:\:j=1,2,\ldots\},$ (\ref{64l}) may be checked explicitly. For the remaining cases we substitute the expression (\ref{64l}) into the recurrence relation  (\ref{62f}) (middle line). The difference between left and right hand sides is	
\begin{equation}\label{64m}	
 \begin{array}{ll}
d(i,i+j+1)+4d(i-1,i+j)-[4i+5]d(i,i+j) = & \\
 -3\times{(-4)^i} \left( s_5(i+1,j)-s_5(i,j)-[4i+5]s_5(i+1,j-1)   \right)    &\\
+2\times {(-4)^{i-1}}  \left( s_9(i,j+1)-s_9(i-1,j+1)-[4i+5]s_9(i,j)   \right),  &  \\
\end{array}
\end{equation}
which is zero on account of (\ref{64aa}), thus concluding the proof.
\end{proof}

\subsection{}

Another interesting and useful problem of a similar nature to that considered in Sec.\:6.1 is discussed in this section. Equation (\ref{21d}) reads
\begin{equation} \label{65a}   
\beta_{n} =
\frac{1}{\Gamma\left(2n+1\right)}\int_{0}^{\infty}\!\!{du}\Phi(u){u^{2n}},\:\:\:n=0,1,\ldots,.
\end{equation}
Suppose we define 
\begin{equation} \label{65b}   
b_{n} =
\int_{0}^{\infty}\!\!{du}\Phi(u){u^{2n}},\:\:\:n=0,1,\ldots,.
\end{equation}
Then from (\ref{21d}) and  (\ref{21g}), with $n>r$, we may write
\begin{equation} \label{65c}				
D(n,r)=\det{\left[\frac{b_{j-i+n}}{\Gamma(2(j-i+n)+1)}\right]}_{i,j=1,\ldots,r}.
\end{equation}
It is therefore of interest to examine the properties of 
\begin{equation} \label{65d}				
\Delta(n,r)=\det{\left[\frac{\Gamma(2n+1)}{\Gamma(2(j-i+n)+1)}\right]}_{i,j=1,\ldots,r}.
\end{equation}

Appendix A of ~\cite{JN4} and some results in ~\cite[p.\,106]{K}    show that $\Delta(n,r)>0,\:\:n>r$. Moreover, define
\begin{equation} \label{65e}   
 \overline{\Delta}(n,r)=\frac{1}{(2n)^{r(r-1)}}{\left( 
{\prod_{i=1}^{r-1} g(i,n)^{r-i}}      \right)}   {\Delta(n,r)},
\end{equation}
where 
\begin{equation} \label{65f}   
g(i,n)=(2n+2i)(2n+2i-1).
\end{equation}
Then we have


\begin{lemma}  \label{L6d}
The function $\overline{\Delta}(n,r)$ is a polynomial in $y=1/(2n)$ of degree $r(r-1)$ for $n>r$. With
\begin{equation} \label{65g}				
\overline{\Delta}(n,r)=\sum_{i=0}^{r(r-1)}\delta(i)y^{i},
\end{equation}
then 
\begin{equation}  \label{65h}
\delta(i)=0,\:\:\:i=0,\ldots,r(r-1)/2-1.
\end{equation}
\end{lemma}
\begin{proof} 
The proof follows a pattern analogous to that used in Lemma  \ref{L6a}.
From the above definitions it follows that
\begin{equation} \label{65i}				
\overline{\Delta}(n,r)_{i,j}={g(j,n-i+1)}\overline{\Delta}(n,r)_{i,j+1},\:\:\:i=1,\ldots,r;\:\:\;
j=1,\ldots,{r-1},
\end{equation}
with $\overline{\Delta}(n,r)_{i,r}=1,\:\:\:i=1,\ldots,r.$

The elements  $ \overline{\Delta}(n,r)_{i,j}$ are polynomials in $y$ of order $2(r-j),\:\:j=1,\ldots,r$. Define a set of column vectors  $\underline{u}(\eta)$ of dimension $r$ by
 $\underline{u}(\eta)=[u_1(\eta),u_2(\eta),\ldots,u_r(\eta)] ^T  $
where 
   \begin{equation}\label{65j}	
   u_{i}(\eta)=(i)^{\eta},\:\:\:i=1,\ldots,r.		
   \end{equation}
Also let the column vector $\underline{U}(j,y)=
 [\overline{\Delta}(n,r)_{1,j},\overline{\Delta}(n,r)_{2,j},\ldots ,  \overline{\Delta}(n,r)_{r,j}] ^T,\:\:\:j=1,\ldots,r,\;\;\:$,
 so that
\begin{equation}\label{65k}	
  \overline{\Delta}(n,r)	  =det\left[  \underline{U}(j,y)
		 \right]	_ {j=1,\ldots,r}
\end{equation}
It may be seen that
\begin{equation}\label{65l}	
    \underline{U}(j,y)=\sum_{k=0}^{2(r-j)}  {c_{j,k}}\: {y^{k}} \underline{u}(k),\:\:\:j=1,\ldots,r
\end{equation}	
for some set of coefficients $c_{j,k}$.

From this point the argument of Lemma  \ref{L6c} applies word for word.
\end{proof}

\section{ $D(n,r)$ FOR LARGE $n$}

\subsection{}

We now turn to the question of the behavior for $D(n,r)$ for large $n$ with $r$ fixed. In ~\cite{JN4} we proposed a method, sketched here, for treating this problem. In Sec.\:7.3 we discuss an alternative method. 

In the current notation it was conjectured that a scaled, normalized version of $D(n,r)$ is approximately
\begin{equation}\label{71a}	
 D(n,r)_{SN} \approx y^{j(r)} \sum_{i=0}^{j(r)} {P(i)}      \left[\frac{1}    {y\beta(n)}\right]^{j(r)-i},\:\:\:n\to {\infty},
\end{equation}
where $j(r)=r(r-1)/2$, $y=(2n)^{-1}$ and $P(i)$ is positive for $i=j(r)$.
 As  $ {n \rightarrow  \infty}$ the quantity   ${1/(y\beta)} {\rightarrow}0$ slowly, so that the term $i=j(r)$ eventually dominates. Our calculations show that other terms can be significant for $n$ in the thousands and beyond, the more so for higher $r$. 

   The key reason for the behavior in (\ref{71a}) is an algebraic result, not special to the case of the RH, which means that certain terms in the expansion ~\cite[(2.26)]{JN4} vanish identically. The relative sizes of the terms in (\ref{71a}) are related to the connection of the formula to the RH. In (\ref{71a}) we have used the Laplace method (see Olver ~\cite{OLV})
to approximate the integrals in (\ref{21d}).

In a later report  ~\cite{JN5}, for the case $r=2$,     we described an improved version of the method of ~\cite{JN4}, which may be generalized to $r>2$ (see Sec.\:7.3). We now present the main steps in the new method for $r=2$. Expanding  (\ref{65c}) gives, for $n>3$ say,
\begin{equation}\label{71b}	
 D(n,2)=  \frac{b_n^2}{\Gamma(2n+1)^2} - \frac{b_{n-1}b_{n+1}}{\Gamma(2n-1)\Gamma(2n+3)}
\end{equation} 
With $g(i,n)$ as defined in (\ref{65f}) 
we normalize $D(n,2)$ to obtain, with $r=2$,
\begin{equation}\label{71d}	
\begin{array} {cll}
 \overline{D}(n,2)  &  = {y^2 } {\left(   {\prod_{i=1}^{r-1} g(i,n)^{r-i} }   \right)}
 { \Gamma(2n+1)^2}   {D(n,r)}    \\     &=  {y^2}    ({b_n^2}(2n+1)(2n+2) - {b_{n-1}b_{n+1}})  \\
  & =B(1)(1+y)(1+2y) - {B(2)}(1-y),
\end{array}
\end{equation}
where we have written $B(1)={b_n^2};\:\:\:B(2)={b_{n-1}b_{n+1}}$.
 We call $B(1),B(2)$ the components of the determinant  $\det{\left[{b_{j-i+n}}\right]}_{i,j=1,\ldots,r}$. 

From  (\ref{65b}), after symmetrizing the double integral (without affecting its value), we may write 
\begin{equation}\label{71e}	
B(k)=\int_{0}^{\infty}\int_{0}^{\infty}dt_1dt_2 \Phi(t_1) \Phi(t_2)t_1^{2n-2}t_1^{2n-2}
R(k,t_1,t_2),\:\:k=1,2,
\end{equation} 
where
\begin{equation}\label{71f}	
R(1,t_1,t_2)= t_1^2t_2^2; \:\:\:       R(2,t_1,t_2)=(t_1^4+t_2^4)/2.
\end{equation} 
Now we set  $t_i=\tau(1+x_i),\:\:i=1,2$, where we will expand the integrand about $t_i= \tau$ to be chosen later. We obtain
\begin{equation}\label{71g}	
\begin{array}{cl}
R(1,t_1,t_2) &   ={\tau^4}[1+2(x_1+x_2)+({x_1}^2+{x_2}^2) +4{x_1}{x_2}  \\
& +  2({x_1}^2{x_2}+{x_2}^2{x_1}) +{x_1}^2{x_2}^2]  \\
  R(2,t_1,t_2) & =({{\tau^4}/2})[2+4(x_1+x_2)+6({x_1}^2+{x_2}^2)   \\
    & + 4({x_1}^3+{x_2}^3) +({x_1}^4+{x_2}^4) ].
\end{array} 
\end{equation} 
The expressions such as $1,(x_1+x_2)/2,({x_1}^2+{x_2}^2)/2,({x_1}{x_2}),
({x_1}^3+{x_2}^3)/2$, etc. in (\ref{71g})
are symmetrized monomials. We can classify them by their type  $ (m,j)$, where $m$ represents the degree and $j=1,\ldots,np(m)$ is the label of a monomial of degree m. For example, for $m=0$,  $np(0)=1$ and the single monomial has structure $1$.  For $m=1$,  $np(1)=1$ and the single monomial has structure $x$. For $m=2$,  $np(2)=2$, and the two types are $x^2$ and $xx$. For $m=3$ we have 2 types $x^3,\:{x^2}x$,   while for $m=4$ there are 3 types $x^4,\: {x^3}x,\:{x^2}{x^2}$. We choose the normalization of each      symmetrized monomial such that, if $x_1=x_2 =x$, then the expression equals $x^m$. We call these monomials $\mu(m,j),\:\:j=1,np(m)$.

With these definitions we  write
\begin{equation}\label{71h}	
R(k,t_1,t_2)={\tau^4}\sum_{m=0}^{4} \sum_{j=1}^{np(m)}T(m,j,k)\mu(m,j),\:\:\:k=1,2,
\end{equation} 
which defines the set of coefficients $\{T(m,j,k)\}$
In view of  (\ref{71e}) and  (\ref{71f}) we may write
\begin{equation}\label{71i}	
 \overline{D}(n,2)=\sum_{i=0}^{2}\sum_{m=0}^{4} \sum_{j=1}^{np(m)}\sum_{k=1}^{2}
     {\epsilon(2,k)}  z(i,k)T(m,j,k){y^i}I(m,j)
\end{equation} 
with
\begin{equation}\label{71j}	
   I(m,j)    =\tau^{4} \int_{0}^{\infty}\int_{0}^{\infty}dt_1dt_2 \Phi(t_1) \Phi(t_2)t_1^{2n-2}t_2^{2n-2}\mu(m,j),
\end{equation} 
and ${\epsilon(2,k)}=1,\:\:k=1;\:\:\:\:=-1,\:\:k=2$.
In  (\ref{71i}) we have used the components of the determinant $ \overline{\Delta}(n,2)$,
\begin{equation}\label{71k}	
Z(1,y)=(1+3y+2y^2):\:\:\:\: Z(2,y)=(1-y),
\end{equation} 
and have written
\begin{equation}\label{71l}	
Z(k,y)=\sum_{i=0}^{2}\sum_{k=1}^{2}z(i,k)y^{i}.
\end{equation} 
Consequently we may write
\begin{equation}\label{71m}	
 \overline{D}(n,2)=\sum_{i=0}^{2}\sum_{m=0}^{4} \sum_{j=1}^{np(m)} C(i;m,j)\:{y^i}I(m,j),
\end{equation} 
where
\begin{equation}\label{71n}	
 C(i;m,j)=\sum_{k=1}^{2}{\epsilon(2,k)} z(i,k)T(m,j,k).
\end{equation} 

To study the behavior  of  $\overline{D}(n,2)$  for large $n$ that is  implied by (\ref{71m}), we note that the coefficients $ \{C(i;m,j)\}$ are independent of $n$, $y=1/(2n)$, and $I(m,j)$ is given by (\ref{71j}). This integral is the product of two integrals each over a single variable of the form
\begin{equation}\label{71o}	
\int_{0}^{\infty}dt \Phi(t) t^{2n-2}{(t-\tau)^{\eta}},
\end{equation} 
with $\eta$ being a small integer. For large $t$   (\ref{21a}) shows that   $\Phi(t)$ is very small, while $t^{2n}$ is large, with the result that the integrand has a sharp maximum at $t=\tau$ where $2{\pi}\exp(4\tau)$. Using the Laplace method (see ~\cite{JN4}, ~\cite{JN5})
we find that the largest terms in (\ref{71m}) occur when $(i,m)=(1,0),\:\:(0,2)$. The ratio of the terms $(1,0)/{(0,2)} =O(\tau), n\to\infty$. Calculations confirm this prediction, including the magnitude of the coefficients of these terms.

It is important to note that $C(0;0,1)=0$, as follows immediately from (\ref{71n}). If this coefficient were not zero, then the ratio $(0,0)/{(1,0)}= O(n),n\to\infty$. To obtain the correct behavior at very large $n$ from (\ref{71n}) it is essential that $C(0;0,1)=0$.

\subsection{}

The method of Sec.\:7.1 for order $2$ can be extended to any $r>2$ in a straightforward manner.   We repeat the main steps for the general case. We write
\begin{equation} \label{72a}
 \overline{D}(n,r)   = {y^{r(r-1)} }   {\left(   \prod_{i=1}^{r-1} {g(i,n)^{r-i} }   \right)}
 { \Gamma(2n+1)^r}   {D(n,r)}  
\end{equation}

Let $q(i,k)),\:\:i=1,2,\ldots,r$ be the permutation of tbe numbers $1,2,\ldots,r$ corresponding to index $k=1,2,\dots,r!$. The order in which these permutations are listed is immaterial.  Define
\begin{equation} \label{72b}				
\nu(i,k)=q(i,k)-i,\:\:i=1,2,\ldots,r.
\end{equation}
 We also define $\epsilon(r,k),\:\:k=1,2,\ldots,r!$ to be equal to $\pm {1}$ depending on whether the permutation is even $(+1)$ or odd $(-1)$.  In place of (\ref{71f}) we choose
\begin{equation} \label{72c}				
R(k,{\underline{t}})={\overline{S}}\prod_{i=1}^{r}{\left[{t_i}^{2(\nu(i,k)+r-1)}\right]},\:\:\:
k=1,2,\ldots,r!
\end{equation}
where ${\underline{t}}=(t_1,t_2,\ldots,t_r)$, and ${\overline{S}}$ means symmetrize.
Then the components of $\det{\left[{b_{j-i+n}}\right]}_{i,j=1,\ldots,r}$ are
   \begin{equation}\label{72d}	
   B(k)=\left\{ \prod_{i=1}^{r}    \int_{0}^{\infty} \!\!  dt_i \Phi(t_i) t_i^{2n-2r+2} \right\}
  R(k,{\underline{t}}),\:\:\:k=1,2,\ldots,r!
  \end{equation} 
For general $r$ the analog of (\ref{71h}) is	
\begin{equation}\label{72e}	
R(k,\underline{t})={\tau^{2r(r-1)}}\sum_{m=0}^{2r(r-1)} \sum_{j=1}^{np(m)}T(m,j,k)\mu(m,j),\:\:\:k=1,2,\ldots,r!
\end{equation} 
which defines the set of coefficients $\{T(m,j,k)\}$ in this case.
The determinant  $\overline{\Delta}(n,r)$ is a polynomial in $y$ of degree $r(r-1)$, as stated  in Lemma  \ref{L6d} above, so for its components we have 
\begin{equation}\label{72f}	
Z(k,y)=\sum_{i=0}^{r(r-1)}\sum_{k=1}^{r!}z(i,k)y^{i}.
\end{equation} 

Combining the above information just as in Sec.\:7.1, we find that
\begin{equation}\label{72g}	
 \overline{D}(n,r)=\sum_{i=0}^{r(r-1)}\sum_{m=0}^{2r(r-1)} \sum_{j=1}^{np(m)} C(i;m,j)\:{y^i}I(m,j),
\end{equation} 
where
\begin{equation}\label{72h}	
 C(i;m,j)=\sum_{k=1}^{r!}{\epsilon(r,k)} z(i,k)T(m,j,k),
\end{equation} 
and
\begin{equation}\label{72i}	
   I(m,j)    =\tau^{2r(r-1)}    \left\{ \prod_{i=1}^{r}    \int_{0}^{\infty} \!\!  dt_i \Phi(t_i) t_i^{2n-2r+2} \right\}       \mu(m,j).
\end{equation} 

Calculations for small $r$ suggest that the behavior noted for $r=2$ may be generalized to all $r$. We propose


\begin{con1}
For $ j=1,\ldots,np(m)$  the coefficients 
\begin{equation} \label{72j}
  C(i;m,j)=0,\:\:m=0,1,\ldots,r(r-1)-2i-2;         \:\:i=0,1,\ldots, r(r-1)/2-1;   
\end{equation}
\end{con1}


\begin{rem3}
The specific values of the cofficients $\{b_n\}$ can have no influence on the proposed relation. Only the structure $b_{j-i+n}$ and the ${\Gamma}$ function in (\ref{63c}	) can matter.
\end{rem3}

We predict that, given the truth of Conjecture 3, the Laplace approach mentioned in Sec.\:7.1 will lead to a rigorous demonstration that,   for fixed order $r$, the determinant $D(n,r)>0$ for sufficiently high $n$.

\subsection{}

We expect that the upper representation will lead to a proof of Conjecture \ref{cnj2}. This has already  been demonstrated\footnote{Will be available at http://publish.uwo.ca/$\sim$jnuttall/upper.pdf} for some of the simpler cases.

Now we return briefly to the approach mentioned at the beginning of Sec.\:6. In   (\ref{24c})  we gave the relation between  $ \Lambda_{[p]}(\underline{s},\underline{t})$
and the compound kernel derived from $K(u,v)=\Phi(u+v)$, i.e.
\begin{equation} \label{66a}				
\Lambda_{[p]}(\underline{s},\underline{t})=
\int_{0}^\infty\!\!\!{d\underline{u}}\int_{0}^\infty\!\!\!{d\underline{v}}\:
\phi_{[p]}  (\underline{u},\underline{s})K_{[p]}(\underline{u},\underline{v})\phi_{[p]}(\underline{v},\underline{t}),
\end{equation}
where \underline{u}, \underline{v}, \underline{s}, \underline{t} are defined as in (\ref{24a}).
In  (\ref{32c}) and  (\ref{33b}) we gave examples showing how $ \Lambda_{[p]}(\underline{s},\underline{t})$ can be related to $D(n,r)$ for $r=3,4$. So long as $n$ is large enough this procedure may be generalized to any value of order $r$. However we know that 
$K(x,y)$ is not sign-regular for $r>4$, so that  (\ref{66a})	cannot be used directly to show that $D(n,r)>0$ for such $r$.

Nevertheless we believe that it may be possible to use  (\ref{66a})	to prove $D(n,r)>0$ so long as $n$ is large enough for a given value of $r$. Just as in our application of the Laplace method we conjecture that, for large $n$, the dominant contribution to  the integral in  (\ref{66a}) arises from locations where all the components of $(\underline{u},\underline{v})$ 
are large. Thus, for large $n$, it should be possible to find a quantity $\mu$ such that  
\begin{equation} \label{66b}				
\Lambda_{[p]}(\underline{s},\underline{t};\mu)=
\int_{\mu}^\infty\!\!\!{d\underline{u}}\int_{\mu}^\infty\!\!\!{d\underline{v}}\:
\phi_{[p]}  (\underline{u},\underline{s})K_{[p]}(\underline{u},\underline{v})\phi_{[p]}(\underline{v},\underline{t}),
\end{equation}
is a good approximation to $ \Lambda_{[p]}(\underline{s},\underline{t})$. Moreover, if $\mu$ is large enough only a small error will be introduced in (\ref{66b}) by replacing $\Phi(u+v)$ in      $K(u,v)$ by the first term in (\ref{21a}). When $\mu$ is large enough that the above approximations are adequate, then we can apply the notions of sign regularity to the reduced region of integration. It follows that $D(n,r)>0$ if Conjecture 3 holds, so that  $\epsilon_r{W_r}(y)>0$ for large enough $n$, where ${W_r}(y)$ is given by  (\ref{61b}).

\section{DISCUSSION}

\subsection{}

The results of this report, both rigorous and conjectured, are a testament to the power of computer-assisted mathematics. In several instances we find that, if we look beneath the surface, there are reasons why the behavior required by the RH occurs. The RH is perhaps not quite as mysterious as it often seems, but no doubt many secrets remain.

 The report raises a number of questions/problems that could led to further progress. We intend to list the problems  at  http://publish.uwo.ca/$\sim$jnuttall, where will also be found the unpublished documents listed in the References.

The author would appreciate comments and corrections.



\bibliographystyle{amsplain}
\bibliography{referenceb}

\end{document}